\documentclass[preprint]{elsarticle}

\usepackage{graphicx}
\usepackage{pdflscape}
\usepackage{amssymb}

\usepackage{lineno}

\usepackage{amsmath}
\usepackage{natbib}

\usepackage{tcolorbox}

\usepackage{caption}  

\usepackage{longtable}
\usepackage[tableposition=t]{caption}
\usepackage{tikz}
\tikzset{near start abs/.style={xshift=1cm}}
\usepackage{multirow}
\usepackage{subfigure}
 \usepackage{diagbox}
\usepackage{amsthm}
\newtheorem{theorem}{Theorem}[section]

\newtheorem{lemma}[theorem]{Lemma}
\theoremstyle{definition}
\newtheorem{definition}{Definition}[section]

\biboptions{authoryear}
\makeatletter
\def\ps@pprintTitle{%
 \let\@oddhead\@empty
 \let\@evenhead\@empty
 \def\@oddfoot{\centerline{\thepage}}%
 \let\@evenfoot\@oddfoot}
\makeatother

\begin{document}

\begin{frontmatter}


\title{An $(R, S)$ Based Heuristic Model for the Stochastic Joint Replenishment Problem}

\author[label1]{Mengyuan Xiang\corref{cor1}}
\ead{mengyuan.xiang@ed.ac.uk}

\author[label1]{Roberto Rossi}
\ead{roberto.rossi@ed.ac.uk}


\author[label2]{S. Armagan Tarim}
\ead{armagan.tarim@ucc.ie}

\cortext[cor1]{Corresponding author}

\address[label1]{Business School, University of Edinburgh, Edinburgh, United Kingdom}
\address[label2]{Cork University Business School, University College Cork, Cork, Ireland}

\begin{abstract}
This paper considers the periodic-review stochastic joint replenishment problem (JRP) under Bookbinder and Tan's static-dynamic uncertainty control policy. According to a static-dynamic uncertainty control rule, the decision maker fixes timing of replenishments once and for all at the beginning of the planning horizon, the inventory position is then raised to a predefined order-up-to-position at the beginning of each replenishment period. In this policy, freezing the replenishment times ameliorates the inherent difficulties pertinent to replenishment coordination of multiple products, whereas dynamic order quantities facilitate dealing with uncertain demands. We adapt and extend an earlier mixed integer linear programming (MILP) model for computing static-dynamic uncertainty policy parameters, and demonstrate that the same can be used to approximate the optimal control rule for the JRP, also known as $(\sigma, \vec{S})$ policy. 
An extensive computational study illustrates the effectiveness of our approach when compared to alternative approaches in the literature. 
\end{abstract}


\end{frontmatter}


\section{Introduction} \label{section1}
The Joint Replenishment Problem (JRP) occurs when several items are ordered from the same supplier, or several products have the same means of transportation, or several products are processed on the same piece of equipment \citep{salamehetal2014}. Every time an order is placed, the group fixed ordering cost is incurred regardless the number of items replenished; in addition there are also item-specific fixed and variable ordering costs that are charged whenever an item is included in a replenishment order. The goal of the JRP is to determine the optimal inventory replenishment plan that minimises the cost of replenishing multiple items. 

The problem of controlling inventory of a multi-item system under joint replenishment has been receiving considerable attention for the past several decades. Literature on JRP can be roughly categorised into deterministic and stochastic based on the nature of demand. In the deterministic joint replenishment inventory system, demand for each individual item is assumed to be constant over an infinite time horizon and replenishments are made at equally spaced time intervals; the problem is to determine the length of replenishment cycles and the frequency of replenishing individual items, e.g., \citep{goyalandBelton1979, kaspiandRosenblatt1991,viswanathan1996,wildemanetal1997,hariga1994,goyal1993,boctoretal2004,nilssonetal2007}. In the stochastic joint replenishment inventory system, the demand for each individual item is unknown but follows certain types of distributions; the problem is to decide the optimal parameters of a given inventory policy, e.g., \citep{balintfy1964, atkinsandIyogun1988, renbergandPlanche1967, kalpakamandArivarignan1993, viswanathan1997, nielsenandLarsen2005, ozkayaetal2006}. Most literature still presents applications to constant and dynamic deterministic demands; however, the study regarding stochastic demand has received increasing attention due to its practical relevance \citep{bastosetal2017}. This work belongs to the growing literature on the stochastic joint replenishment.

This paper applies the static-dynamic strategy, proposed by \cite{bt1988} for tackling single-item lot-sizing problems, in the context of a JRP system. The static-dynamic strategy, known as $(R, S)$, features two control parameters: $R$, timing of replenishment, and $S$, order-up-to-position. At each review period, the decision maker places an order so as to increase the inventory position (i.e. net inventory level + outstanding orders) to a given order-up-to-position. In the context of the JRP system, a periodic-review $(R, S)$ policy is adopted for each item. The $(R, S)$ policy is an appealing strategy since it eases the coordination between supply chain players \citep{kilicandtarim2011}, and facilitates managing joint replenishment \citep{silveretal1998}. Additionally, the $(R, S)$ policy comes with the advantage of being able to tackle nonstationary demand which has not been addressed yet in the literature.

Our goal is to tackle the periodic-review nonstationary JRP under $(R, S)$ policy. We first present a mixed-integer linear programming (MILP) model for computing policy parameters that minimise the expected total cost comprising group fixed ordering costs, item-specific fixed ordering costs, holding costs, and penalty costs over the planning horizon. Our model generalises the discussion in \citep{rossietal2015}, which presented an MILP model for approximating optimal $(R, S)$ policy parameters for single-item lot-sizing problems. 
We further show that our MILP model can be used to approximate the $(\sigma, \vec{S})$ policy, which is known to be optimal for the JRP \citep{liuandesogbue2012}. Under this policy, decision makers order up to $\vec{S}$ if opening inventory positions fall in $\sigma$ ($\sigma \subset R^N, \vec{S} \in R^N$, $N$ represents the number of items) at the beginning of each time period. 
Numerical experiments illustrate the effectiveness of our policy and the corresponding MILP model. 

We contribute to the literature on the stochastic JPR as follows.
\begin{itemize}
\item We present, for the first time in the literature, a mathematical programming (MP) model for tackling the {\em nonstationary} stochastic JRP.
\item We reformulate this MP model as an MILP model that can be solved using off-the-shelf solvers. 
\item We demonstrate that our MILP model can be used to approximate the $(\sigma, \vec{S})$ policy, which is optimal for the JRP.
\item In an extensive computational study based on existing test beds, drawn from the literature, we demonstrate the effectiveness of our models when compared to alternative approaches in the literature.
\end{itemize}

This rest of this paper is organised as follows. Section \ref{literaturereview} surveys relevant literature. Section \ref{problemdescription} describes problem settings. Section \ref{milp} presents an MILP model for computing $(R, S)$ policy parameters. Section \ref{kconvexity} extends the MILP model for approximating the optimal $(\sigma, \vec{S})$ policy parameters. An extensive computational study is conducted in Section \ref{computationalstudy}. We draw conclusions in Section \ref{conclusion}. 

\section{Literature review}\label{literaturereview}
The problem of controlling the inventory of a multi-item system under joint replenishment has received increasing attention over the past several decades. For a thorough review of literature readers may refer to \citep{silverandPeterson1985,goyalandSatir1989,vanetal1992, khoujaandGoyal2008, bastosetal2017}. In this section, we focus our attention on existing policies for tackling stochastic JRPs. In particular, we survey control policies that have been considered in the literature.

{\bf $(\sigma, \vec{S})$ policy}. The landmark study of \cite{Scarf1960} proved the optimality of ($s$,$S$) policies for the single-item inventory problem; since then, there have been attempts to generalise this result to multi-item inventory systems. \cite{johnson1967} characterised the optimal policy for the stationary case and introduced the $(\sigma, \vec{S})$ policy, in which $\sigma \subset \mathcal{R}^N$ and $\vec{S} \in \mathcal{R}^N$; in this policy one orders up to $\vec{S}$ if inventory levels $\vec{I} \in \sigma$ and $\vec{I} \leq \vec{S}$, otherwise one does not order. \cite{kalin1980} showed that, when $\vec{I} \in \sigma$ and $\vec{I} \nleq \vec{S}$, there exists $\vec{S}(\vec{I}) \geq \vec{I}$ such that the optimal policy is to order up to $\vec{S}(\vec{I})$, this policy is named $(\sigma, \vec{S}(\cdot))$ policy. \cite{ohnoetal1994} proposed an algorithm for computing an optimal ordering policy $(\sigma,\vec{S}(\cdot))$ for a periodic-view multi-item inventory system.  \cite{ohnoandishigaki2001} further proposed a policy iteration method to compute an exact optimal policy by leaving properties of the optimal policy for continuous-time inventory problems with compound Poisson demands. 
\cite{gallegoandsethi2005} gave the general definition of $K$-convexity in $\mathcal{R}^N$, which encompasses both the joint ordering and individual ordering case; it derived an optimal policy for the two-item deterministic inventory problem with a joint ordering cost. However, the computation of the optimal $(\sigma, \vec{S})$ policies is still a difficult task.



{\bf $(s, c, S)$ policy.} Several works on stochastic JRPs have focused on computing $(s, c, S)$ policies, introduced by \cite{balintfy1964}. This policy features three control parameters: $s$, reorder point; $c$, can-order level; $S$, order-up-to-position. Under this policy, When the inventory position of item $n$ crosses $s_n$, a replenishment order is triggered to raise its inventory position to $S_n$; meanwhile,  any other item $j$ with an inventory position at or below its can-order point, $c_j (s_j < c_j < S_j)$, is also included in the replenishment, raising its inventory position to $S_j$. Under the assumption of Poisson-distributed demands, \cite{ignall1969} proved that the $(s, c, S)$ policy is not optimal even for two-item problems. \cite{silver1974} proposed the decomposition method to approximate $(s, c, S)$ policy parameters, where the multi-item problem is decomposed into several single-item problems. This approximation technique was followed by \citep{melchiors2002, johansenandMelchiors2003}. \cite{kayisetal2008} modelled the two-item JRP problem as a semi-Markov decision model, and proposed an enumerative approach to approximate $(s, c, S)$ policies. In addition, \citep{schaackandSilver1972, thompstoneandSilver1975, silver1981, federgruenetal1984} studied JRPs with compound Poisson-distributed demands.



{\bf $(R, T)$ policy.}
\cite{atkinsandIyogun1988} proposed two periodic-review $(R, T)$-type policies, namely periodic policy $P$ and modified periodic policy $MP$, which differ only in the way the ordering periods $T_n$ are determined. Under this policy, every $T_n$ periods, the inventory position of item $n$ is raised to $R_n$. Numerical experiments demonstrate that the $MP$ policy performs consistently better than the $(s, c, S)$ policy, and that the $P$ policy generally outperforms the the $(s, c, S)$ policy, except for problems involving small values of group fixed ordering cost.

{\bf $(Q, S)$ policy.} This policy was first proposed by \cite{renbergandPlanche1967}. Under this policy,  whenever the total inventory position drops to the group reorder point, an order is placed to raise inventory position of each item to item-specific order-up-to-position $S$. The combined order quantity is $Q$, and the group reorder point is reached when the combined usage reaches $Q$. \cite{pantumsinchai1992} evaluated the computational performance of the $(Q, S)$ policy by comparing it against the $(s, c, S)$ policy, $P$ policy and  $MP$ policy on the basis of long-run total average costs. Computational experiments showed that the $MP$ policy consistently outperforms the $(s, c, S)$ policy on the test instances, and both $MP$ and $(Q, S)$ policy perform better as the group ordering cost increases. The study showed that the $(Q, S)$ policy is appropriate for items for which the stock-out costs are low and the major set-up cost is high relative to the minor set-up cost.


{\bf $P(s, S)$ policy.} This policy was proposed by \cite{viswanathan1997} for periodic-review inventory systems, in which inventory position of each item is reviewed at every fixed and constant time interval. At each review time, the $(s, S)$ policy is applied to each item, so that any item with inventory position at or below $s$ is order up to $S$. For a fixed review period, the algorithm of \cite{zhengandFedergruen1991} is adopted to compute the optimal $(s, S)$ policy parameters. Computational studies indicated that although the proposed policy requires more computational effort, it generally dominates the $MP$ policy, and dominates $(s, c, S)$ policy, and $(Q, S)$ policy for most test instances.

{\bf $Q(s, S)$ policy.}  \cite{nielsenandLarsen2005} combined features of $(Q, S)$ policy and $P(s, S)$ policy, and proposed the $Q(s, S)$ policy. By operating under this policy, the total inventory position is continuously reviewed while the item-specific inventory positions are reviewed only when the total consumption since the last order reaches $Q$. Then every item with inventory position less than or equal to its respective reorder point $s$ is order to $S$. An analytical solution is derived by using the Markov decision theory in \cite{nielsenandLarsen2005}. Computational study demonstrated that the $Q(s, S)$ policy outperforms $P(s, S$) policy, and dominates $(Q, S)$ policy in $17$ of $18$ test instances on the data set of \cite{atkinsandIyogun1988}.  

{\bf $(Q, S, T)$ policy.} This continuous-review policy was proposed by \cite{ozkayaetal2006}. Decision makers raise the inventory position of each item $i$ to its order-up-to-position $S_i$ whenever a total of $Q$ demands accumulated or $T$ time units have elapsed, whichever occurs first. This policy is a hybrid of the continuous review $(Q, S)$ policy, proposed by \cite{renbergandPlanche1967}, and the periodic review $(R, T)$ policy, proposed by \cite{atkinsandIyogun1988}. Thus, it features benefits of two separate policies. The comprehensive numerical study indicates that the proposed policy dominates the $P(s, S)$ policy, $(Q, s)$ policy, $Q(s, S)$ policy, and $(s, c, S)$ policy in $100$ of $139$ instances. 

{\bf $(R, S)$ policy} was proposed by \cite{bt1988} for controlling single-item inventory systems. This policy requires decision makers to place an order at each replenishment period to increase the inventory position to the order-up-to-position $S$. It is an appealing policy since it eases the coordination between supply chain players \citep{kilicandtarim2011}, and facilitates managing joint replenishment \citep{silveretal1998}. \citep{tarimandkingsman2004, tarimandkingsman2006} formulated a mixed integer programming (MIP) model for computing optimal $(R, S)$ policy parameters. \cite{tarimetal2011} relaxed the  MIP model, and solved it as a shortest path problem which does not require the use of a mathematical programming solver. In addition, \cite{ozenetal2012} introduced a DP-based algorithm for solving small-size problems, and an approximation heuristic and a relaxation heuristic for tackling larger-size problems. \cite{rossietal2015} generalised the discussions above and developed a unified MILP model for approximating the $(R, S)$ policy by adopting the piecewise linear approximation technique in \cite{rossietal2014}. Recently, \cite{tuncetal2018} presented an extended MIP model that blends heuristic methods originally introduced by \cite{tuncetal2014} and \cite{rossietal2015}. As a result, this formulation features the computational efficiency of \cite{tuncetal2014} and the modelling variety of \cite{rossietal2015}. Although various efficient modelling methods for computing $(R, S)$ policy parameters have been proposed, all existing works focus on a single-item inventory system. 


The stochastic JRP is an open research area for the development of more efficient computational methods and control policies. The main purpose of this work is to apply the $(R, S)$ policy to a multi-item inventory system. In the context of the JRP, we apply the periodic-review $(R,S)$ policy, originally proposed by \cite{bt1988} for tackling single-item lot sizing problems, to JRPs under stochastic demand and fixed lead time; a periodic-review $(R, S)$ policy is adopted for each item. Note that when the stochastic demand is stationary, the $(R, S)$ policy is the same as the $MP$ policy proposed by \cite{atkinsandIyogun1988}, in which every $T_n$ periods one raises the inventory position of item $n$ to the order-up-to-position $R_n$. However, the $(R, S)$ policy can also deal with non-stationary stochastic demands; a setting that was not addressed in \cite{atkinsandIyogun1988}. 
In what follows, we introduce a novel MILP approach for approximating $(R, S)$ policies under non-stationary stochastic demands for multi-item inventory systems. Nonlinear costs are approximated by leveraging the technique originally introduced in \citep{rossietal2014}.

\section{Problem description}\label{problemdescription}
Consider a periodic-review $N$-item inventory management system over a $T$-period planning horizon. We assume that demand $d_t^n$ of item $n$, $n=1, \ldots, N$, in period $t$, $t=1, \ldots, T$ is a random variable with a known probability density function; all $d_t^n$ are assumed to be mutually independent.

We further assume that replenishments are issued at the beginning of each time period. There is a group fixed ordering cost $K$, which is incurred whenever a replenishment is issued regardless the number of items replenished. Moreover, there is an item-specific fixed ordering cost $k^n$, which is incurred whenever item $n$ is replenished regardless the quantity of the replenishment. 

We define $Q_t^n$ as the quantity of item $n$ ordered in period $t$, which is placed and received immediately. Then, the ordering cost of item $n$ in period $t$ with ordering quantity $Q_t^n$ can be written as,
\begin{align}
&c_t^n(Q_t^n)=\begin{cases} k^n, &Q_t^n > 0,\\
0, &Q_t^n=0.\end{cases}
\end{align}
Let $c_t(\vec{Q}_t)$ denote the ordering cost of period $t$ with ordering quantity vector $\vec{Q}_t = (Q_t^1, \ldots, Q_t^N)$. $c_t(\vec{Q}_t)$ has the following structure
\begin{align}
&c_t(\vec{Q}_t)=\begin{cases}K+ \sum_{n=1}^N c_t^n(Q_t^n), &\exists Q_t^n|Q_t^n >0,\\
0,& \mbox{otherwise}. \end{cases}
\end{align}

A penalty cost $b^n$ is incurred for each unit of item $n$ of backorder demand per period, and a holding cost $h^n$ is charged for each unit of item $n$ carried from one period to the next. The immediate penalty and holding cost of period $t$ can be expressed as
\begin{align}
&L_t(\vec{y})=\sum_{t=1}^n\Big(b^n\cdot \text{E}[\max(d_t^n-y^n, 0)]+h^n\cdot \text{E}[\max(y^n-d_t^n, 0)]\Big),
\end{align}
where vector $\vec{y}=(y^1, \ldots, y^N)$ is the inventory level immediately after replenishments are received at the beginning of period $t$, and ``$\text{E}$" denotes the expectation operator.

Let $I_t^n$ denote the net inventory level of item $n$ at the end of period $t$, which is also the opening inventory level of period $t+1$, and $C_t(\vec{I}_{t-1})$ denote the expected total cost of the optimal plan over period $t, \ldots, T$, given opening inventory level $\vec{I}_{t-1}=(I_{t-1}^1, \ldots, I_{t-1}^N)$ at the beginning of period $t$. Then, $C_t(\vec{I}_{t-1})$ can be written as, 
\begin{align}
\small
C_t(\vec{I}_{t-1})=
\min_{\vec{Q}_t}\big\{ c_t(\vec{Q}_t)+ L_t(\vec{I}_{t-1}+\vec{Q}_t)+E[C_{t+1}(\vec{I}_{t-1}+\vec{Q}_t-\vec{D}_t)]\big\}
\end{align}
where $\vec{D}_t=(d_t^1, \cdots, d_t^N)$, and
\begin{align} 
C_T(\vec{I}_{T-1})=
\min_{\vec{Q}_T} \big\{c_T(\vec{Q}_T)+ L_T(\vec{I}_{T-1}+\vec{Q}_{T})\big\},
\end{align}
represents the boundary condition.

\begin{figure}[!htbp]
\begin{center}
\includegraphics[width=8.4cm]{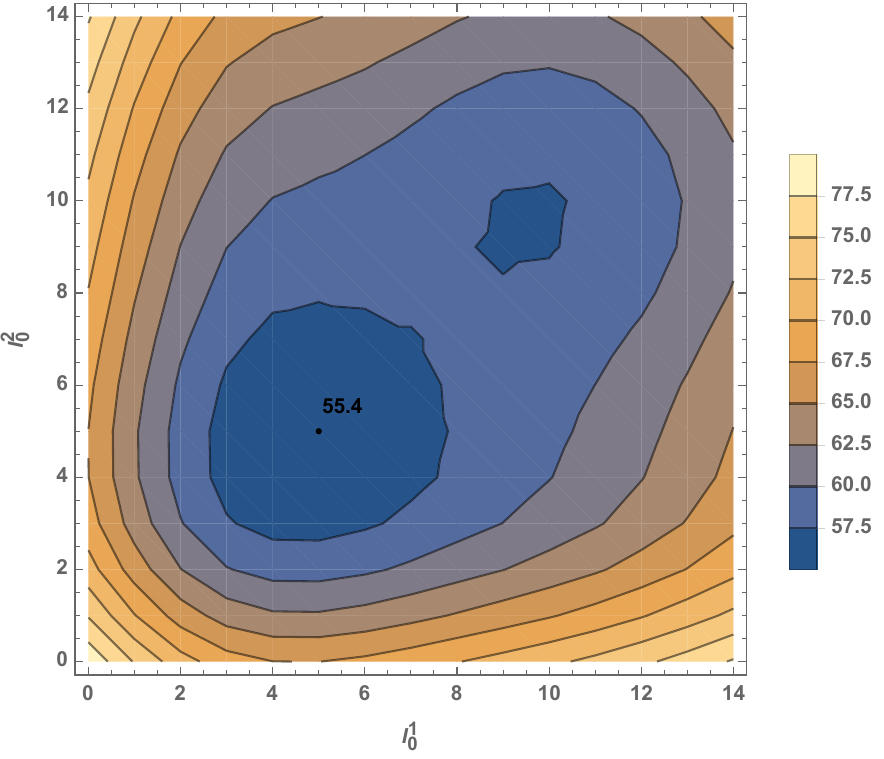}    
\caption{Expected total cost, i.e. $G_1(\vec{I}_0)$, contour plot for the two-item joint replenishment numerical example} 
\label{fig:example_countour_plot}
\end{center}
\end{figure}

{\bf Example.} We consider an instance with two items in which the group fixed ordering cost is $K=10$, and for each item, the item-specific ordering cost $k$ is 0, the holding cost is $h=1$, and the stock-out penalty cost is $b=5$. We control the inventory for two items over a planning horizon of $T=4$ periods. We assume that the demand of item $n$ in period $t$ follows a Poisson distribution with rate $\lambda^n_t$; where $\lambda^1_t=\lambda^2_t=\{3,6,9,6\}$. The expected total cost, i.e. $C_1(\vec{I}_0)$, of an optimal policy, given initial inventory level $I_0^1=I_0^2=0$, can be obtained via stochastic dynamic programming (SDP) and is equal to $65.4$. In Fig. \ref{fig:example_countour_plot} we plot $G_1(\vec{I}_0)$ for $I_0^1 \in [0, 14]$ and $I_0^2 \in [0, 14]$.



\section{An MILP model for approximating non-stationary stochastic $(R, S)$ policies}\label{milp}
In this section, we formulate the stochastic JRP problem under the $(R, S)$ policy as an MILP model. Under the $(R, S)$ policy, the replenishment periods and associated order-up-to-positions are fixed at the beginning of the planning horizon, while actual order quantities are decided at the beginning of each replenishment period. Note that in the context of the JRP, a periodic-review $(R, S)$ policy is adopted for each item. We next introduce a stochastic programming formulation (Section \ref{stochasticprogramming}) and then approximate it by an MILP model (Section \ref{MINLPformulation}).

\subsection{A stochastic program}
\label{stochasticprogramming}
Consider the periodic-review $N$-item $T$-period JRP described in Section \ref{problemdescription}. We introduce binary variables $\delta_t$ and $y_t^n$, $t=1, \ldots, T$, and $n=1, \ldots, N$; $\delta_t$ takes value $1$ if a group order is made in period $t$, otherwise $0$; $y_t^n$ is set to $1$ if item $n$ is replenished in period $t$. 

Under the $(R, S)$ policy, we reformulate the stochastic dynamic programming model in Section \ref{problemdescription} as the stochastic program in Fig. \ref{spformulation}. 
\begin{figure}[!htbp]
\small
\begin{align}
&\min
\sum_{t=1}^T\Big(K\cdot\delta_t+\sum_{n=1}^N(k^n\cdot y_t^n+b^n\text{E}[\max(-I_t^n, 0)]+h^n\text{E}[\max(I_t^n, 0)])\Big)\label{sp}
\end{align}
Subject to, $n=1, \ldots, N$, $t=1, \ldots, T$, 
\begin{align}
&\delta_t \geq y_t^n  \label{sp-1}\\
&I_t^n = I_0^n +\sum_{i=1}^{t}Q_i^n-\sum_{j=1}^td_j^n  \label{sp-2}\\
&y_t^n =\begin{cases}
1, &Q_t^n > 0,\\
0, &Q_t^n = 0.
\end{cases}\label{sp-5}\\
&Q_t^n \geq 0   \label{sp-3}\\
&\delta_t =\{0, 1\}  \label{sp-4}\\
&I_t^n \in \mathcal{R}  \label{sp-6}
\end{align}
\caption{Stochastic programming formulation for the JRP.}
\label{spformulation}
\end{figure}
 
The objective is to find the optimal replenishment plan so as to minimise the expected ordering costs, penalty costs, and holding costs of $N$ items over the $T$-period planning horizon. Constraints (\ref{sp-1}) imply that if at least an item is ordered, then a group replenishment is issued. Constraints (\ref{sp-2}) are inventory conservation constraints: the inventory level at the end of period $t$ is equal to the initial inventory level, plus all orders received before the end of period $t$, minus demands raised up to period $t$.  Constraints (\ref{sp-5})- (\ref{sp-6}) state domains of $y_t^n$, $Q_t^n$, $\delta_t$, and $I_t^n$. 
 
 \subsection{MILP model for approximating $(R, S)$ policies}\label{MINLPformulation}
By leveraging the piecewise approximation approach in \citep{rossietal2014}, the stochastic programming formulation in Fig. \ref{spformulation} can be approximated by an MILP model. 

We introduce the first order loss function  \[\mathcal{L}(x,\omega)=\text{E}[\max(\omega-x,0)],\] and its complementary function \[\hat{\mathcal{L}}(x,\omega)=\text{E}[\max (x-\omega,0)],\] where $\omega$ is a random variable with a known probability density function, and $x$ is a scalar variable.

Consider a single replenishment cycle of item $n$ over periods $i, \ldots, j$, where the only replenishment is placed at the beginning of period $i$ with order-up-to-level $S_i^n$, and the initial inventory level is $I_{i-1}^n$. Thus, $I_t^n$, $t=i, \ldots, j$, must equal to the order-up-to-level $S_{i}^n$, minus the demand convolution $d^n_{i,t}$ over periods $i, \ldots, t$, i.e.: $I_t^n=S_{i}^n-d_{i,t}^n$. We rewrite the expected excess back-orders $\max(-I_t^n, 0)$ and on-hand stocks $\max(I_t^n, 0)$ as $\mathcal{L}(S_i^n, d_{i,t}^n)$ and $\hat{\mathcal{L}}(S_i^n, d_{i,t}^n)$, by means of the first order loss function and its complementary function.

We introduce a binary variable $P_{jt}^n$, $j=1, \ldots, t$, $t=1, \ldots, T$, $n=1, \ldots, N$, which is set to one if the most recent replenishment of item $n$ up to period $t$ was issued in period $j$, where $j\leq t$ --- if no replenishment occurs before or at period $t$, then we let $P_{1t}^n=1$, this allows us to properly account for demand variance from the beginning of the planning horizon. We observe that if $P_{jt}^n=1$, the closing inventory level of period $t$ must equal to the order-up-to-level of period $j$, minus the demand convolution over periods $j, \ldots, t$, i.e. $I_t^n=S_j^n-d_{jt}^n$. Then, the expected excess back-orders and on-hand stocks of period $t$ can be written by means of the first order loss function and its complementary function, $\sum_{j=1}^t\mathcal{L}(S_j^n,d_{jt}^n)P_{jt}^n$, and $\sum_{j=1}^t\hat{\mathcal{L}}(S_j^n, d_{jt}^n)P_{jt}^n$. Additionally, since period $j$ must be the only most recent order received up to period $t$, the following constraints must be satisfied.
 \begin{align}
 &\sum_{j=1}^tP_{jt}^n=1,& \label{P-1}\\
 &P_{jt}^n\geq y_j^n-\sum_{k=j+1}^t y_k^n, & j=1, \ldots, t. \label{P-2}
 \end{align}

In what follows, let ``$\sim$" denote the expectation operator. We introduce decision variables $\tilde{B}_t^n \geq 0$ and $\tilde{H}_t^n\geq 0$ to represent expected excess back-orders and on-hand stocks. The stochastic program in Fig \ref{spformulation} can be approximated by the MINLP model in Fig. \ref{MINLPmodel}. 
\begin{figure}[!htbp]
\small
\begin{align}
    \min \sum_{t=1}^T\Big(K\cdot \delta_{t}+\sum_{n=1}^N\big(k^n\cdot y_t^n+h^n\tilde{H}_{t}^n+b^n\tilde{B}_{t}^n\big)\Big) \label{MINLP-0}
\end{align}
Subject to, $n = 1, \ldots, N$, $t=1, \ldots, T$,
\begin{align}
&\delta_t\geq y_t^n & 
\label{MINLP-1}\\
&\tilde{I}_t^n+\tilde{d}_t^n-\tilde{I}_{t-1}^n\geq 0 &
\label{MINLP-2}\\
&y_t^n=0 \rightarrow \tilde{I}_t^n+\tilde{d}_t^n-\tilde{I}_{t-1}^n= 0 
\label{MINLP-3}\\
&\sum_{j=1}^tP_{jt}^n=1 & 
\label{MINLP-4}\\
&P_{j,t}^n \geq y_{j}^n - \sum_{k=j+1}^{t}y_k^n & j={1, \ldots, t} \label{MINLP-5}\\
&P_{jt}^n=1 \rightarrow \tilde{H}_t^n=\hat{L}(\tilde{I}_t^n+\tilde{d}_{j,t}^n,d_{j,t}^n) & j={1, \ldots, t} \label{MINLP-6}\\
&P_{jt}^n=1 \rightarrow 
\tilde{B}_t^n=L(\tilde{I}_t^n+\tilde{d}_{j,t}^n,d_{j,t}^n)& j={1, \ldots, t}
\label{MINLP-7}\\
&\delta_t =\{0, 1\} & 
\label{MINLP-8}\\
&y_t^n =\{0, 1\} & 
\label{MINLP-9}\\
&P_{jt}^n=\{0, 1\} & j={1, \ldots, t}&
\label{MINLP-10}
\end{align}
\caption{MINLP model for approximating $(R, S)$ policies}
\label{MINLPmodel}
\end{figure}

The objective function (\ref{MINLP-0}) minimizes the expected group fixed ordering costs, item-specific fixed ordering costs, holding costs, and penalty costs of $N$-item over the $T$-period planning horizon. Constraints (\ref{MINLP-1}) imply an individual item can only be included in a group replenishment if that replenishment is made. Constraints (\ref{MINLP-2}) - (\ref{MINLP-3}) are inventory balance constraints. Constraints (\ref{MINLP-4}) - (\ref{MINLP-5}) ensure the most recent order before period $t$ was issued in period $j$. Nonlinear constraints (\ref{MINLP-6}) - (\ref{MINLP-7}) represent the expected on-hand stocks and back-orders of item $n$ over the planning horizon. Note that the order-up-to-position of item $n$ in period $j$ can be expressed by the expected closing inventory level and expected demand convolution, i.e., $S_{j}^n = \tilde{I}_t^n+\tilde{d}_{j,t}^n$. Constraints (\ref{MINLP-8}) - (\ref{MINLP-10}) indicate domains of $\delta_t$, $y_t^n$, and $P_{jt}^n$.

By solving the model in Fig. \ref{MINLPmodel}, the optimal replenishment plan including group replenishment periods $\delta_t$, and item-specific replenishment periods $y_t^n$, and the item-specific order-up-to-positions $S_{t}^n=\tilde{I}_{t}^n+\tilde{d}_{t}^n$ are obtained, for $t=1, \ldots, T$, and $n=1, \ldots, N$. The MINLP model in Fig. \ref{MINLPmodel} can be readily approximated by an MILP model by using the approach discussed in \citep{rossietal2014,rossietal2015} to piecewise linearise loss functions in constraints (\ref{MINLP-6}) and (\ref{MINLP-7}).
Note that the MINLP model in Fig. \ref{MINLPmodel} can be extended to explore the fixed lead time settings. Details are discussed in \ref{leadtime}

{\bf Example.} We demonstrate the modelling strategy underpinning the MILP model on a $5$-item $10$-period example. It is assumed that demands follow a Poisson distribution with rates $\lambda_t^n$ presented in Table \ref{demandrate}. The initial inventory level is taken as zero,  and item-specific lead time $L^n=[1, 2, 3,  1, 3]$. Other parameters are: $K=500$, $b=10$, $h=2$, and $k^n=[120, 100, 80, 120, 150]$.  We employ eleven segments in the piecewise-linear approximations of $B_{t}^n$ and $H_{t}^n$ (for $n=1,\ldots, 5$, and $t=1,\ldots, 10$).
\begin{table}[!ht]
    \centering
    \begin{tabular}{|c|cccccccccc|}
    \hline
    \diagbox{item}{$\lambda_t^n$}{period}& 1     & 2     & 3     & 4     & 5     & 6     & 7     & 8     & 9     & 10 \\
    \hline
    1     & 40    & 40    & 40    & 40    & 40    & 40    & 40    & 40    & 40    & 40 \\
    2     & 5     & 64    & 29    & 54    & 70    & 50    & 54    & 45    & 13    & 50 \\
    3     & 40    & 55    & 72    & 86    & 78    & 51    & 42    & 38    & 30    & 26 \\
    4     & 41    & 58    & 75    & 63    & 40    & 35    & 33    & 18    & 29    & 39 \\
    5     & 45    & 40    & 22    & 31    & 38    & 46    & 59    & 62    & 46    & 40 \\
    \hline
    \end{tabular}
    \caption{Demand rates $\lambda_t^n$ of the $5$-item $10$-period example}
    \label{demandrate}
\end{table}

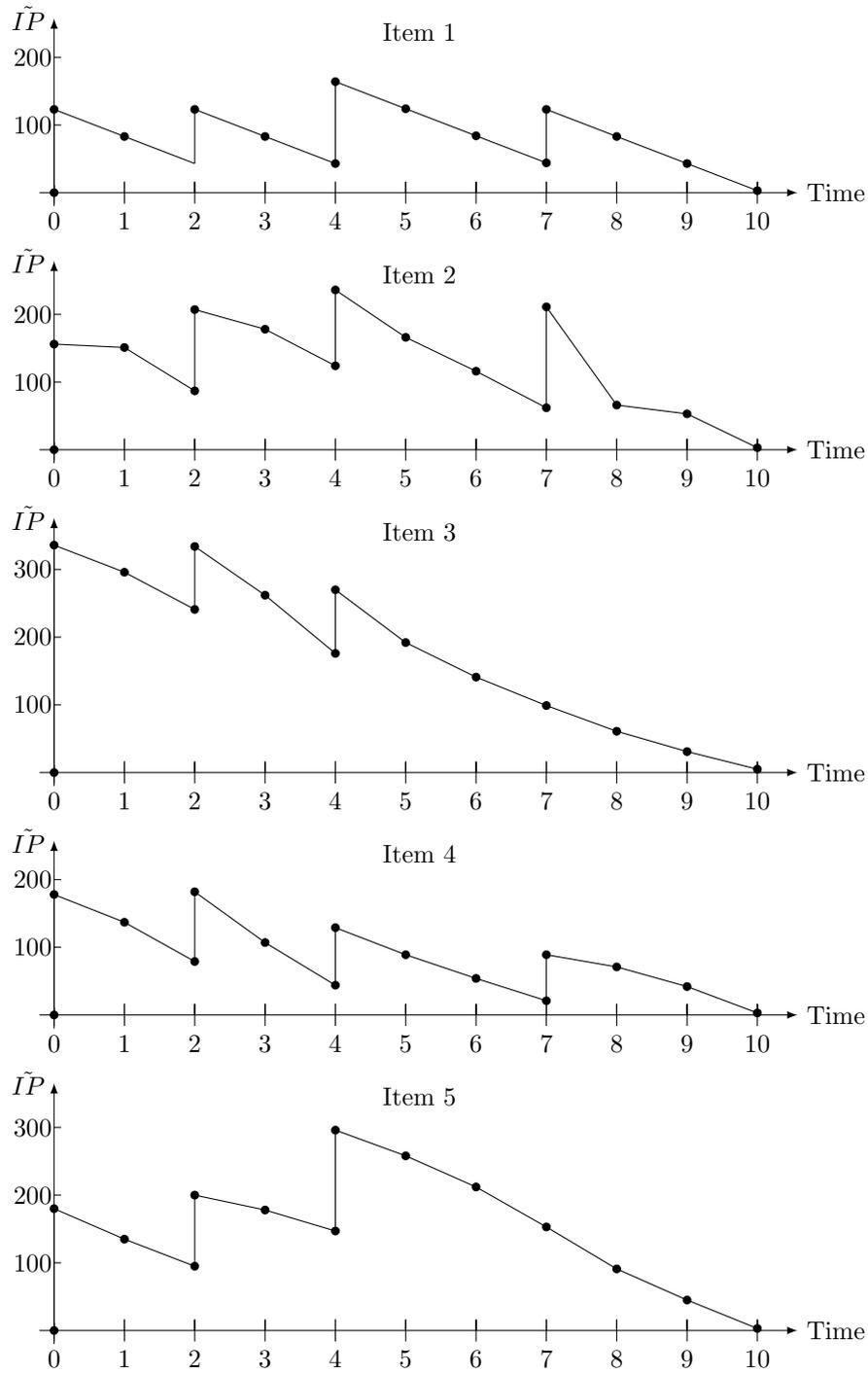
\begin{figure}[!htbp]
\centering
\begin{tikzpicture}[x=0.0480952380952381cm, y=0.018518518518518517cm]
 \draw [-latex] ([xshift=-2mm] 0.0,0) -- ([xshift=100mm] 3.5,0) node[right] {Time};
 \draw (0.0,0) -- +(0mm,1mm) -- +(0mm,-1.5mm) node[below] {0};
 \draw (20,0) -- +(0mm,1.5mm) -- +(0mm,-1.5mm) node[below] {1};
 \draw (40,0) -- +(0mm,1.5mm) -- +(0mm,-1.5mm) node[below] {2};
 \draw (60,0) -- +(0mm,1.5mm) -- +(0mm,-1.5mm) node[below] {3};
 \draw (80,0) -- +(0mm,1.5mm) -- +(0mm,-1.5mm) node[below] {4};
 \draw (100,0) -- +(0mm,1.5mm) -- +(0mm,-1.5mm) node[below] {5};
 \draw (120,0) -- +(0mm,1.5mm) -- +(0mm,-1.5mm) node[below] {6};
 \draw (140,0) -- +(0mm,1.5mm) -- +(0mm,-1.5mm) node[below] {7};
 \draw (160,0) -- +(0mm,1.5mm) -- +(0mm,-1.5mm) node[below] {8};
 \draw (180,0) -- +(0mm,1.5mm) -- +(0mm,-1.5mm) node[below] {9};
 \draw (200,0) -- +(0mm,1.5mm) -- +(0mm,-1.5mm) node[below] {10};
 \draw [-latex] ([yshift=-0mm] 0,-10.0) -- ([yshift=22mm] 0, 10.0) node[left] {$\tilde{IP}$};
 
\draw (0,50) -- (2,50) node[left] {100};
\draw (0,100) -- (2,100) node[left] {200};

\filldraw[black] (0,0) circle (1.5pt);
\filldraw[black] (0,61.5) circle (1.5pt);
\filldraw[black] (20,41.5) circle (1.5pt);
\filldraw[black] (40,61.5) circle (1.5pt);
\filldraw[black] (60,41.5) circle (1.5pt);
\filldraw[black] (80,21.5) circle (1.5pt);
\filldraw[black] (80,82) circle (1.5pt);
\filldraw[black] (100,62) circle (1.5pt);
\filldraw[black] (120,42) circle (1.5pt);
\filldraw[black] (140,22) circle (1.5pt);
\filldraw[black] (140,61.5) circle (1.5pt);
\filldraw[black] (160,41.5) circle (1.5pt);
\filldraw[black] (180,21.5) circle (1.5pt);
\filldraw[black] (200,1.5) circle (1.5pt);

 \draw (0,0)--(0,61.5);
 \draw (0,61.5)--(20,41.5);
\draw (20,41.5)--(40,21.5);
\draw (40,21.5)--(40,61.5);   
\draw (40,61.5)--(60,41.5); 
\draw (60,41.5)--(80,21.5);
\draw (80,21.5)--(80,82);  
\draw (80,82)--(100,62);  
\draw (100,62)--(120,42); 
\draw (120,42)--(140,22); 
\draw (140,22)--(140,61.5); 
\draw (140,61.5)--(160,41.5); 
\draw (160,41.5)--(180,21.5); 
\draw (180,21.5)--(200,1.5); 

\node at (50mm,22mm) {Item 1};

\end{tikzpicture}
\\
\begin{tikzpicture}[x=0.0480952380952381cm, y=0.018518518518518517cm]
 \draw [-latex] ([xshift=-2mm] 0.0,0) -- ([xshift=100mm] 3.5,0) node[right] {Time};
 \draw (0.0,0) -- +(0mm,1mm) -- +(0mm,-1.5mm) node[below] {0};
 \draw (20,0) -- +(0mm,1.5mm) -- +(0mm,-1.5mm) node[below] {1};
 \draw (40,0) -- +(0mm,1.5mm) -- +(0mm,-1.5mm) node[below] {2};
 \draw (60,0) -- +(0mm,1.5mm) -- +(0mm,-1.5mm) node[below] {3};
 \draw (80,0) -- +(0mm,1.5mm) -- +(0mm,-1.5mm) node[below] {4};
 \draw (100,0) -- +(0mm,1.5mm) -- +(0mm,-1.5mm) node[below] {5};
 \draw (120,0) -- +(0mm,1.5mm) -- +(0mm,-1.5mm) node[below] {6};
 \draw (140,0) -- +(0mm,1.5mm) -- +(0mm,-1.5mm) node[below] {7};
 \draw (160,0) -- +(0mm,1.5mm) -- +(0mm,-1.5mm) node[below] {8};
 \draw (180,0) -- +(0mm,1.5mm) -- +(0mm,-1.5mm) node[below] {9};
 \draw (200,0) -- +(0mm,1.5mm) -- +(0mm,-1.5mm) node[below] {10};
 \draw [-latex] ([yshift=-0mm] 0,-10.0) -- ([yshift=24mm] 0, 10.0) node[left] {$\tilde{IP}$};
 
\draw (0,50) -- (2,50) node[left] {100};
\draw (0,100) -- (2,100) node[left] {200};

\filldraw (0,0) circle (1.5pt);
\filldraw(0,78) circle (1.5pt);
\filldraw (20,75.5) circle (1.5pt);
\filldraw (40,43.5) circle (1.5pt);
\filldraw (40,103.5) circle (1.5pt);
\filldraw (60,89) circle (1.5pt);
\filldraw (80,62) circle (1.5pt);
\filldraw (80,118) circle (1.5pt);
\filldraw (100,83) circle (1.5pt);
\filldraw (120,58) circle (1.5pt);
\filldraw (140,31) circle (1.5pt);
\filldraw (140,105.5) circle (1.5pt);
\filldraw (160,33) circle (1.5pt);
\filldraw (180,26.5) circle (1.5pt);
\filldraw (200,1.5) circle (1.5pt);

\draw (0,0)--(0,78);
\draw (0,78)--(20,75.5);
\draw (20,75.5)--(40,43.5);
\draw (40,43.5)--(40,103.5);   
\draw (40,103.5)--(60,89); 
\draw (60,89)--(80,62);
\draw (80,62)--(80,118);  
\draw (80,118)--(100,83);  
\draw (100,83)--(120,58); 
\draw (120,58)--(140,31); 
\draw (140,31)--(140,105.5); 
\draw (140,105.5)--(160,33); 
\draw (160,33)--(180,26.5); 
\draw (180,26.5)--(200,1.5); 

\node at (50mm,24mm) {Item 2};

\end{tikzpicture}
\\
\begin{tikzpicture}[x=0.0480952380952381cm, y=0.018518518518518517cm]
 \draw [-latex] ([xshift=-2mm] 0.0,0) -- ([xshift=100mm] 3.5,0) node[right] {Time};
 \draw (0.0,0) -- +(0mm,1mm) -- +(0mm,-1.5mm) node[below] {0};
 \draw (20,0) -- +(0mm,1.5mm) -- +(0mm,-1.5mm) node[below] {1};
 \draw (40,0) -- +(0mm,1.5mm) -- +(0mm,-1.5mm) node[below] {2};
 \draw (60,0) -- +(0mm,1.5mm) -- +(0mm,-1.5mm) node[below] {3};
 \draw (80,0) -- +(0mm,1.5mm) -- +(0mm,-1.5mm) node[below] {4};
 \draw (100,0) -- +(0mm,1.5mm) -- +(0mm,-1.5mm) node[below] {5};
 \draw (120,0) -- +(0mm,1.5mm) -- +(0mm,-1.5mm) node[below] {6};
 \draw (140,0) -- +(0mm,1.5mm) -- +(0mm,-1.5mm) node[below] {7};
 \draw (160,0) -- +(0mm,1.5mm) -- +(0mm,-1.5mm) node[below] {8};
 \draw (180,0) -- +(0mm,1.5mm) -- +(0mm,-1.5mm) node[below] {9};
 \draw (200,0) -- +(0mm,1.5mm) -- +(0mm,-1.5mm) node[below] {10};
 \draw [-latex] ([yshift=-0mm] 0,-10.0) -- ([yshift=33mm] 0, 10.0) node[left] {$\tilde{IP}$};
 
\draw (0,50) -- (2,50) node[left] {100};
\draw (0,100) -- (2,100) node[left] {200};
\draw (0,150) -- (2,150) node[left] {300};

\filldraw[black] (0,0) circle (1.5pt);
\filldraw[black] (0,168) circle (1.5pt);
\filldraw[black] (20,148) circle (1.5pt);
\filldraw[black] (40,120.5) circle (1.5pt);
\filldraw[black] (40,167) circle (1.5pt);
\filldraw[black] (60,131) circle (1.5pt);
\filldraw[black] (80,88) circle (1.5pt);
\filldraw[black] (80,135) circle (1.5pt);
\filldraw[black] (100,96) circle (1.5pt);
\filldraw[black] (120,70.5) circle (1.5pt);
\filldraw[black] (140,49.5) circle (1.5pt);
\filldraw[black] (160,30.5) circle (1.5pt);
\filldraw[black] (180,15.5) circle (1.5pt);
\filldraw[black] (200,2.5) circle (1.5pt);

\draw (0,0)--(0,168);
 \draw (0,168)--(20,148);
\draw (20,148)--(40,120.5);
\draw (40,120.5)--(40,167);   
\draw (40,167)--(60,131); 
\draw (60,131)--(80,88);
\draw (80,88)--(80,135);  
\draw (80,135)--(100,96);  
\draw (100,96)--(120,70.5); 
\draw (120,70.5)--(140,49.5); 
\draw (140,49.5)--(160,30.5); 
\draw (160,30.5)--(180,15.5); 
\draw (180,15.5)--(200,2.5); 

\node at (50mm,33mm) {Item 3};

\end{tikzpicture}
\\
\begin{tikzpicture}[x=0.0480952380952381cm, y=0.018518518518518517cm]
  \draw [-latex] ([xshift=-2mm] 0.0,0) -- ([xshift=100mm] 3.5,0) node[right] {Time};
 \draw (0.0,0) -- +(0mm,1mm) -- +(0mm,-1.5mm) node[below] {0};
 \draw (20,0) -- +(0mm,1.5mm) -- +(0mm,-1.5mm) node[below] {1};
 \draw (40,0) -- +(0mm,1.5mm) -- +(0mm,-1.5mm) node[below] {2};
 \draw (60,0) -- +(0mm,1.5mm) -- +(0mm,-1.5mm) node[below] {3};
 \draw (80,0) -- +(0mm,1.5mm) -- +(0mm,-1.5mm) node[below] {4};
 \draw (100,0) -- +(0mm,1.5mm) -- +(0mm,-1.5mm) node[below] {5};
 \draw (120,0) -- +(0mm,1.5mm) -- +(0mm,-1.5mm) node[below] {6};
 \draw (140,0) -- +(0mm,1.5mm) -- +(0mm,-1.5mm) node[below] {7};
 \draw (160,0) -- +(0mm,1.5mm) -- +(0mm,-1.5mm) node[below] {8};
 \draw (180,0) -- +(0mm,1.5mm) -- +(0mm,-1.5mm) node[below] {9};
 \draw (200,0) -- +(0mm,1.5mm) -- +(0mm,-1.5mm) node[below] {10};
 \draw [-latex] ([yshift=-0mm] 0,-10.0) -- ([yshift=22mm] 0, 10.0) node[left] {$\tilde{IP}$};
 
\draw (0,50) -- (2,50) node[left] {100};
\draw (0,100) -- (2,100) node[left] {200};

\filldraw[black] (0,0) circle (1.5pt);
\filldraw[black] (0,89) circle (1.5pt);
\filldraw[black] (20,68.5) circle (1.5pt);
\filldraw[black] (40,39.5) circle (1.5pt);
\filldraw[black] (40,91) circle (1.5pt);
\filldraw[black] (60,53.5) circle (1.5pt);
\filldraw[black] (80,22) circle (1.5pt);
\filldraw[black] (80,64.5) circle (1.5pt);
\filldraw[black] (100,44.5) circle (1.5pt);
\filldraw[black] (120,27) circle (1.5pt);
\filldraw[black] (140,10.5) circle (1.5pt);
\filldraw[black] (140,44.5) circle (1.5pt);
\filldraw[black] (160,35.5) circle (1.5pt);
\filldraw[black] (180,21) circle (1.5pt);
\filldraw[black] (200,1.5) circle (1.5pt);

\draw (0,0)--(0,89);
 \draw (0,89)--(20,68.5);
\draw (20,68.5)--(40,39.5);
\draw (40,39.5)--(40,91);   
\draw (40,91)--(60,53.5); 
\draw (60,53.5)--(80,22);
\draw (80,22)--(80,64.5);  
\draw (80,64.5)--(100,44.5);  
\draw (100,44.5)--(120,27); 
\draw (120,27)--(140,10.5); 
\draw (140,10.5)--(140,44.5); 
\draw (140,44.5)--(160,35.5); 
\draw (160,35.5)--(180,21); 
\draw (180,21)--(200,1.5); 

\node at (50mm,22mm) {Item 4};

\end{tikzpicture}
\\
\begin{tikzpicture}[x=0.0480952380952381cm, y=0.018518518518518517cm]
 \draw [-latex] ([xshift=-2mm] 0.0,0) -- ([xshift=100mm] 3.5,0) node[right] {Time};
 \draw (0.0,0) -- +(0mm,1mm) -- +(0mm,-1.5mm) node[below] {0};
 \draw (20,0) -- +(0mm,1.5mm) -- +(0mm,-1.5mm) node[below] {1};
 \draw (40,0) -- +(0mm,1.5mm) -- +(0mm,-1.5mm) node[below] {2};
 \draw (60,0) -- +(0mm,1.5mm) -- +(0mm,-1.5mm) node[below] {3};
 \draw (80,0) -- +(0mm,1.5mm) -- +(0mm,-1.5mm) node[below] {4};
 \draw (100,0) -- +(0mm,1.5mm) -- +(0mm,-1.5mm) node[below] {5};
 \draw (120,0) -- +(0mm,1.5mm) -- +(0mm,-1.5mm) node[below] {6};
 \draw (140,0) -- +(0mm,1.5mm) -- +(0mm,-1.5mm) node[below] {7};
 \draw (160,0) -- +(0mm,1.5mm) -- +(0mm,-1.5mm) node[below] {8};
 \draw (180,0) -- +(0mm,1.5mm) -- +(0mm,-1.5mm) node[below] {9};
 \draw (200,0) -- +(0mm,1.5mm) -- +(0mm,-1.5mm) node[below] {10};
 \draw [-latex] ([yshift=-0mm] 0,-10.0) -- ([yshift=32mm] 0, 10.0) node[left] {$\tilde{IP}$};
 
\draw (0,50) -- (2,50) node[left] {100};
\draw (0,100) -- (2,100) node[left] {200};
\draw (0,150) -- (2,150) node[left] {300};

\filldraw[black] (0,0) circle (1.5pt);
\filldraw[black] (0,90) circle (1.5pt);
\filldraw[black] (20,67.5) circle (1.5pt);
\filldraw[black] (40,47.5) circle (1.5pt);
\filldraw[black] (40,100) circle (1.5pt);
\filldraw[black] (60,89) circle (1.5pt);
\filldraw[black] (80,73.5) circle (1.5pt);
\filldraw[black] (80,148) circle (1.5pt);
\filldraw[black] (100,129) circle (1.5pt);
\filldraw[black] (120,106) circle (1.5pt);
\filldraw[black] (140,76.5) circle (1.5pt);
\filldraw[black] (160,45.5) circle (1.5pt);
\filldraw[black] (180,22.5) circle (1.5pt);
\filldraw[black] (200,1.5) circle (1.5pt);

\draw (0,0)--(0,90);
\draw (0,90)--(20,67.5);
\draw (20,67.5)--(40,47.5);
\draw (40,47.5)--(40,100);   
\draw (40,100)--(60,89); 
\draw (60,89)--(80,73.5);
\draw (80,73.5)--(80,148);  
\draw (80,148)--(100,129);  
\draw (100,129)--(120,106); 
\draw (120,106)--(140,76.5); 
\draw (140,76.5)--(160,45.5); 
\draw (160,45.5)--(180,22.5); 
\draw (180,22.5)--(200,1.5); 

\node at (50mm,32mm) {Item 5};

\end{tikzpicture}
\\
\caption{Replenish plans of the $5$-item $10$-period example}
\label{replenishplans}
\end{figure}

The resulting expected total cost is $14236$. Replenishment plans of each item are presented in Fig. \ref{replenishplans}. Items $1$, $2$ and $4$ are replenished in periods $1$, $3$, $5$, and $8$; while item $3$ and $5$ are replenished only in periods $1$, $3$, and $5$, since orders in period $8$ could not be received by the end of the planning horizon. Additionally, since item $1$ is expected to be ordered every two periods with the same order-up-to-position $123$ by the nature of stationary demand, while it is ordered up to a higher position $164$ in period $5$ to cover demands in the next $3$ periods in order to coordinate with other items. 

\section{MILP model for approximating the optimal ($\sigma, \vec{S}$) policies} \label{kconvexity}
Since the landmark study of \cite{Scarf1960} which proved the optimality for the single-item inventory system, there have been several attempts to prove the optimality for multi-item inventory systems, e.g.: \citep{johnson1967, kalin1980, ohnoandishigaki2001, gallegoandsethi2005}. In this section we demonstrate that the MILP model proposed in Section \ref{MINLPformulation} can be used to approximate the optimal replenishment plan under $(\sigma, \vec{S})$ policy for the JRP. 

 \begin{definition}
 Function $f(\cdot): \mathcal{R}^N\rightarrow \mathcal{R}$ is $K$-convex if 
\[f(ax+(1-a)z) \leq a f(x)+(1-a)[f(z)+{K}\delta(z-x)],\]
where $x \leq z$, $a \in [0,1]$, and ${K}\delta(z-x)$ is defined as follows,
\[{K}\delta(z-x)=K\delta(e'x)+\sum_{n=1}^Nk^n\delta(x_n),\]
where $e'=(1, 1, \cdots, 1)'\in \mathcal{R}^N$, $\delta(0)=0$, and $\delta(y)=1$ for all $y>0$.
 \end{definition}
 
\vspace{1em} 

\cite{gallegoandsethi2005} showed the optimal policy for the joint setup cost case by studying the function 
\begin{align}
 G_t(\vec{y}) =L_t(\vec{y})+ C_{t+1}(\vec{y}-\vec{d}_t). 
\end{align}

Consider a continuous $K$-convex function $G_t(\cdot)$, 
then it has global minimum at $\vec{S}_t$. Define set $\Sigma = \{\vec{I}_{t-1}\leq \vec{S}_t| G_t(\vec{I}_{t-1}) \leq G_t(\vec{S}_t)+K\}$, and set $\sigma=\{\vec{I}_{t-1}\leq \vec{S}_t|\vec{I}_{t-1} \notin \Sigma\}$. Lemma \ref{lemma_K-convexity} shows that the optimal replenishment plan is to order up to $\vec{S}_t$ if opening inventory levels $\vec{I}_{t-1} \in \sigma$ and $\vec{I}_{t-1} \leq \vec{S}_t$; and not to order, otherwise.

\vspace{1em} 
\begin{lemma}[\cite{gallegoandsethi2005}]\label{lemma_K-convexity}
If $G$ is continuous $K$-convex, continuous and coercive, then
\begin{itemize}
    \item $\vec{I} \in \Sigma \Rightarrow G(\vec{I}) \leq K+G(\vec{S})$,
    \item $\vec{I} \in \sigma \Rightarrow G(\vec{I}) > K+G(\vec{S})$.
\end{itemize}
\end{lemma}


We next show that the MILP model in Fig. \ref{MINLPmodel} can be adjusted to approximate set $\sigma$ and $\vec{S}$. 

Due to the complexity of $\sigma$, it is impractical to derive a closed form expression for it. To address this difficulty, {\em we propose a strategy to determine whether a given initial inventory level vector} $\vec{I}_0$ {\em belongs to} $\sigma$. By solving our modified MILP model over the planning horizon $k, \ldots, T$, we observe the minimised expected total cost $G_k(\vec{S}_k)$, order-up-to-levels $\vec{S}_k$, and the first period replenishment decision $\delta_k$. If $\delta_k =1$, then $\vec{I}_{k-1} \in \sigma$; otherwise, $\vec{I}_{k-1} \in \Sigma$. Therefore, our MILP model can be used to determine whether given initial inventory levels $\vec{I}_0 \in \sigma$. Moreover, by repeating this procedure, one can approximate the optimal replenishment strategy for every period $k=1, \ldots, T$.

\begin{figure}[!ht]
\centering
\subfigure[Expected total cost contour plot obtained via MILP approximation] {
\label{fig:example_countour_plot_milp}
\includegraphics[width=0.45\textwidth]{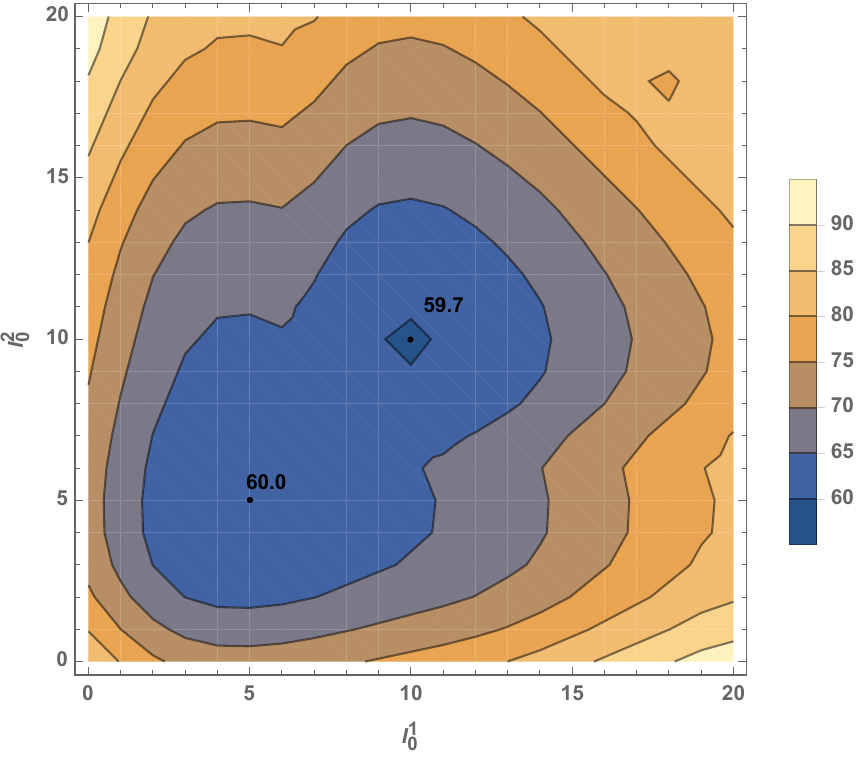}
}
\subfigure[Plot of expected total costs obtained via MILP and SDP]{
\label{fig:example_optimality}
\includegraphics[width=0.45\textwidth]{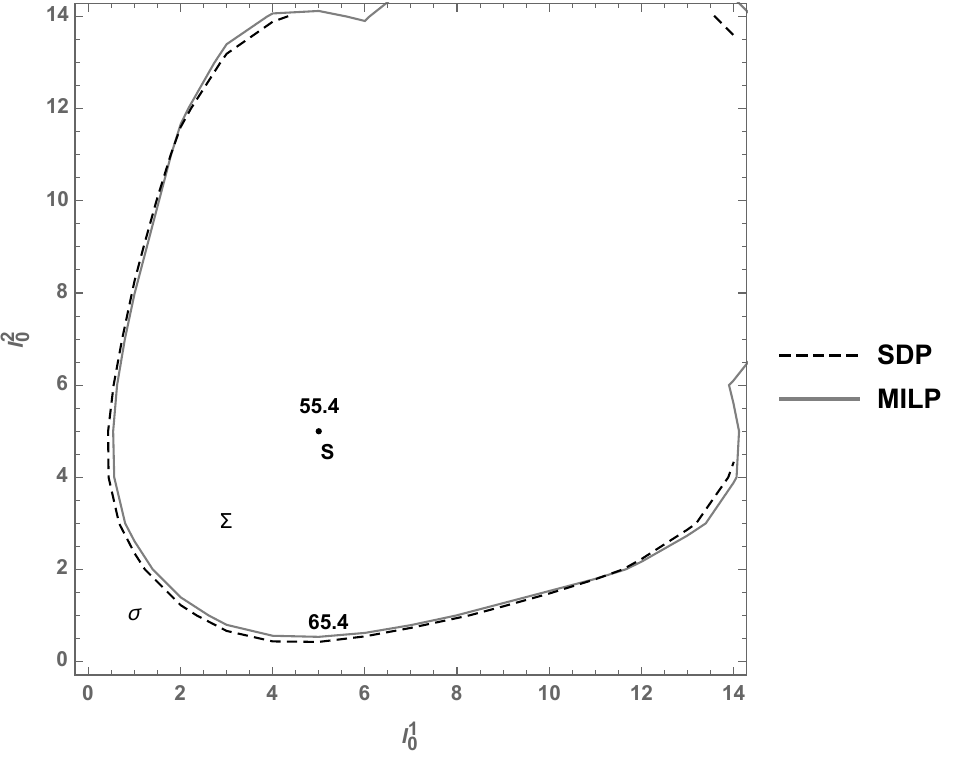}
}
\caption{Plot of expected total costs for the two-item joint replenishment numerical example}
\end{figure}

{\bf Example.} We illustrate the concept introduced on the $2$-item $4$-period example presented in Section \ref{problemdescription}. Assuming initial inventory levels $\vec{I}_0 \in [0, 1, \ldots, 20]\times [0, 1, \ldots, 20]$, we plot the expected total cost contours, obtained via the modified MILP, in Fig. \ref{fig:example_countour_plot_milp}. Note that there are two similar minima, which is expected, since the ordering cost is relatively small and the demand variance is large. We plot set $\sigma$ and $\vec{S}$ obtained via the modified MILP model, and compare them with that obtained via the stochastic dynamic programming in Fig. \ref{fig:example_optimality}. The optimal replenishment plan is to place an order whenever inventory levels $\vec{I}_0=(I_0^1, I_0^2)$ fall in set $\sigma$, and not to place an order if $\vec{I}_0$ fall in $\Sigma$. We observe that set $\sigma$ and $\vec{S}$ obtained via the modified MILP model neatly approximate those obtained via the stochastic dynamic programming.


\section{Computational Experiments}\label{computationalstudy}
In this section we assess the cost performance of the $(R, S)$ policy by comparing its cost performance against $(Q, S, T)$ policy \citep{ozkayaetal2006}, $Q(s, S)$ policy \citep{nielsenandLarsen2005}, $P(s, S)$ policy \citep{viswanathan1997}, $(Q, S)$ policy \citep{pantumsinchai1992},  $MP$ policy \citep{atkinsandIyogun1988}, $(s, c, S)_M$ policy \citep{melchiors2002}, and $(s, c, S)_F$ policy \citep{federgruenetal1984}, on data sets of \cite{atkinsandIyogun1988} and \cite{viswanathan1997}. These data sets consider stationary demand over an infinite horizon. Unfortunately, computing $(R, S)$ policy parameters for infinite horizon JRPs via our MILP model is computationally expensive; however, since the demand is stationary, it is possible to derive an efficient shortest path reformulation, which we present in \ref{shortespathapproximation} and we use in our computational study. 

Computational experiments are conducted by using IBM ILOG CPLEX Optimization Studio 12.7 and Matlab R2016a on a 3.20 GHz Intel Core i5-6500 CPU with 16.0 GB RAM, 64 bit machine. 

Since the shortest path reformulation operates over a finite horizon, in order to compare the cost performance of the $(R, S)$ policy with the continuous-review $(s, c, S)$, $(Q, S)$, and $(Q, S, T)$ policies, we discretize each time period into $20$ small periods. We consider a planning horizon length of $6.6$ periods for a total of $132$ small periods. For each test instance, we first obtain the optimal replenishment plan by solving the shortest path reformulation presented in  \ref{shortespathapproximation}. The computational time is limited to $5$ minutes, if a timeout occurs, the best solution available is adopted. Next, we simulate the expected average cost of each test instance via Monte Carlo Simulation (100,000 replications). Finally, we compare the average cost per small period  against the average cost under existing policies.

The data set of \cite{atkinsandIyogun1988} assumes that the demand of each item follows a stationary Poisson distribution with rate $\lambda^n$, $n=1, \ldots, 12$. The item-specific fixed ordering cost $k^n$, expected demand $\lambda^n$, and lead time $L^n$ are displayed in Table \ref{set1parameters}. Items share the same penalty cost $b=30$, holding cost $h \in \{2, 6, 20\}$, and group fixed ordering cost $K \in \{20, 50, 100, 150, 500\}$.

\begin{table}[!htbp]
    \centering
    \begin{tabular}{c|c|c|c|c|c|c|c|c|c|c|c|c}
    \hline
         items& 1 & 2& 3& 4& 5& 6& 7& 8& 9& 10& 11& 12  \\
          \hline
         $k^n$& 10& 10& 20& 20& 40& 20& 40& 40& 60& 60& 80& 80 \\
         $\lambda^n$& 40& 35& 40& 40& 40& 20& 20& 20& 28& 20& 20& 20\\
         $L^n$&0.2& 0.5& 0.2& 0.1& 0.2& 1.5& 1.0& 1.0& 1.0& 1.0& 1.0& 1.0\\
         \hline
    \end{tabular}
    \caption{$K^n$, $\lambda^n$, and $L^n$ of data set \cite{atkinsandIyogun1988}}
    \label{set1parameters}
\end{table}
 
The data set of \cite{atkinsandIyogun1988} contains some unusual lot sizing instances; more specifically, instances for which the group as well as item fixed ordering costs become negligible in comparison to holding costs. In the lot-sizing literature the fixed ordering cost is commonly assumed to be greater than the holding cost \citep[see][p. 62, Property 2]{citeulike:8526547}; moreover, the penalty cost should not be smaller than the holding cost. Additionally, the fixed ordering cost should be greater than the penalty cost, otherwise, the inventory system tends to place orders in every period instead of penalising backorders. To focus on meaningful lot sizing instances --- instances in which a trade off between fixed ordering and holding/penalty cost is sought --- we filter test instances of the data set of \cite{atkinsandIyogun1988} by using the following conditions: $K > b \geq h$. We also check the order frequency in each period and we discard instances in which orders are issued too frequently --- i.e. instance in which a replenishment is issued more than twice per time period, as it turns out that for these instances order coordination is straightforward due to negligible item fixed ordering costs: if a group order is placed, all items are ordered. We present the filtered computational results in Table \ref{dataset1}. 
\begin{table*}[!htbp]
  \centering
  \resizebox{0.95\linewidth}{!}{
       \begin{tabular}{lll|l|rrrrrrr}
    \hline
    \multicolumn{1}{r}{\multirow{2}[4]{*}{$K$}} & \multicolumn{1}{r}{\multirow{2}[4]{*}{$b$}} & \multicolumn{1}{r|}{\multirow{2}[4]{*}{$h$}} & \multicolumn{1}{r|}{\multirow{2}[4]{*}{$(R, S)$}} & \multicolumn{7}{c}{Average cost improvement $\Delta\%$} \\
\cline{5-11}         
&       & \multicolumn{1}{r|}{} &       & $(Q, S, T)$ & $Q(s, S)$ & $P(s, S)$ & $(Q, S)$ & $MP$ & $(s, c, S)\_M$ & $(s, c, S)\_F$ \\
    \hline
    \multicolumn{1}{r}{50} & \multicolumn{1}{r}{30} & \multicolumn{1}{r|}{2} & \multicolumn{1}{r|}{936.94} & \textcolor[rgb]{ 1,  0,  0}{\textbf{-0.91}} & \textcolor[rgb]{ 1,  0,  0}{\textbf{-0.84}} & \textcolor[rgb]{ 1,  0,  0}{\textbf{-0.33}} & 4.38  & 0.68  & 0.79  & 2.14  \\
    \multicolumn{1}{r}{100} & \multicolumn{1}{r}{30} & \multicolumn{1}{r|}{2} & \multicolumn{1}{r|}{990.50} & \textcolor[rgb]{ 1,  0,  0}{\textbf{-0.05}} & \textcolor[rgb]{ 1,  0,  0}{\textbf{-0.45}} & 0.75  & 2.57  & 1.77  & 4.39  & 6.81  \\
    \multicolumn{1}{r}{150} & \multicolumn{1}{r}{30} & \multicolumn{1}{r|}{2} & \multicolumn{1}{r|}{1046.56} & \textcolor[rgb]{ 1,  0,  0}{\textbf{-0.24}} & \textcolor[rgb]{ 1,  0,  0}{\textbf{-1.01}} & \textcolor[rgb]{ 1,  0,  0}{\textbf{-0.35}} & 0.52  & 0.65  & 5.68  & 8.36  \\
    \multicolumn{1}{r}{200} & \multicolumn{1}{r}{30} & \multicolumn{1}{r|}{2} & \multicolumn{1}{r|}{1072.97} & 1.32  & 0.47  & 1.11  & 1.34  & 2.12  & 8.34  & 12.31  \\
    \multicolumn{1}{r}{100} & \multicolumn{1}{r}{30} & \multicolumn{1}{r|}{6} & \multicolumn{1}{r|}{1639.75} & \textcolor[rgb]{ 1,  0,  0}{\textbf{-0.23}} & \textcolor[rgb]{ 1,  0,  0}{\textbf{-1.52}} & \textcolor[rgb]{ 1,  0,  0}{\textbf{-1.02}} & 2.15  & 0.00  & 1.24  & 3.31  \\
    \multicolumn{1}{r}{150} & \multicolumn{1}{r}{30} & \multicolumn{1}{r|}{6} & \multicolumn{1}{r|}{1707.05} & 0.64  & \textcolor[rgb]{ 1,  0,  0}{\textbf{-0.60}} & \textcolor[rgb]{ 1,  0,  0}{\textbf{-0.07}} & 1.46  & 0.95  & 2.34  & 6.68  \\
    \multicolumn{1}{r}{200} & \multicolumn{1}{r}{30} & \multicolumn{1}{r|}{6} & \multicolumn{1}{r|}{1766.38} & 1.16  & 0.08  & 0.65  & 1.17  & 1.67  & 3.08  & 9.04  \\
    \multicolumn{1}{r}{150} & \multicolumn{1}{r}{30} & \multicolumn{1}{r|}{20} & \multicolumn{1}{r|}{2718.47} & 0.77  & 4.32  & \textcolor[rgb]{ 1,  0,  0}{\textbf{-1.26}} & 1.27  & \textcolor[rgb]{ 1,  0,  0}{\textbf{-0.21}} & \textcolor[rgb]{ 1,  0,  0}{\textbf{-0.59}} & 6.20  \\
    \multicolumn{1}{r}{200} & \multicolumn{1}{r}{30} & \multicolumn{1}{r|}{20} & \multicolumn{1}{r|}{2812.52} & \textcolor[rgb]{ 1,  0,  0}{\textbf{-3.23}} & 0.14  & \textcolor[rgb]{ 1,  0,  0}{\textbf{-0.72}} & 0.77  & 0.34  & 0.25  & 8.34  \\
    \hline
    \multicolumn{4}{l|}{Average cost improvement  $\Delta\%$} & \textcolor[rgb]{ 1,  0,  0}{\textbf{-0.09}} & 0.07  & \textcolor[rgb]{ 1,  0,  0}{\textbf{-0.14}} & 1.74  & 0.89  & 2.84  & 7.02  \\
  \hline
    \end{tabular}%
  }
    \caption{Computational results on the data set of \cite{atkinsandIyogun1988}}
  \label{dataset1}%
\end{table*}
 
 Let $\Delta \%$ denote the percentage gap between the expected average cost of existing policies and that of the proposed $(R, S)$ policy, over the expected average cost of the $(R, S)$ policy. By definition, a positive $\Delta\%$ represents the $(R, S)$ policy outperforms existing policies.  
 Note that expected average costs under $(Q, S, T)$, $Q(s, S)$, $P(s, S)$, $(Q, S)$, and $(s,c, S)_M$ policies are obtained from \cite{ozkayaetal2006}, that of $(s, c, S)_F$ policy is obtained from \cite{melchiors2002}, and that of $MP$ policy is obtained from \cite{viswanathan1997}. 
 
 We observe that the $(R, S)$ policy fully dominates all policies in $2$ of $9$ test instances; $(Q, S, T)$ is the best policy in $2$ instances; $Q(s, S)$ is the best policy in $4$ instances; $P(s, S)$ is the best policy in 1 instance. Moreover, the $(R, S)$ policy outperforms $(Q, S)$ and $(s, c, S)_F$ policies, and no policy is dominant on all test instances. The average cost improvement $\Delta\%$ increases with the increase of group fixed ordering cost, and decreases with the increase of holding cost compared with $(s, c, S)_M$ and $(s, c, S)_F$ policies. That means an increase in group fixed ordering cost or a decrease in holding cost improves the cost performance of $(R, S)$ policy. If we compare the $(R, S)$ policy with $(Q, S, T)$, $Q(s, S)$, $P(s, S)$, $(Q, S)$, and $MP$ policies, there is no obvious trend with respect to the group fixed ordering cost and holding cost. The $(R, S)$ policy performs better than $Q(s, S)$, $(Q, S)$, $MP$, $(s, c, S)_M$, and $(s, c, S)_F$ policies with average improvements of $0.07\%$, $1.74\%$, $0.89\%$, $2.84\%$, and $7.02\%$, respectively; however, $(Q, S, T)$ and $P(s, S)$ policies perform slightly better than the $(R, S)$ policy with average improvements of $0.09\%$ and $0.14\%$, respectively.
 
\cite{viswanathan1997} adopts the same experimental setup as \cite{atkinsandIyogun1988}, except $h \in \{2, 6, 10, 200, 600, 1000\}$, $K \in \{20, 50, 100, 200, 500\}$, and $b \in \{10, 50, 100, 200, 1000, 5000, 10000, 20000\}$. 

We filter the computational results by using the same conditions previously adopted. We present computational results of the $(R, S)$ policy on the data set of \cite{viswanathan1997} in Table \ref{dataset2}. 
We observe that the $(R, S)$ policy dominates $13$ of $31$ test instances; $(Q, S, T)$ is the best policy in $13$ instances; $Q(s, S)$ is the best policy in $9$ instances; $P(s, S)$ is the best policy in $1$ instances. There is once more no dominant policy on all test instances. Regarding the comparison with other policies, the average cost improvement $\Delta\%$ decreases as the penalty cost increases; while there is no obvious trend with respect to the group fixed ordering cost, and penalty cost. On average, the $(R, S)$ policy performs better than $Q(s, S)$, $P(s, S)$, $(Q, S)$, $MP$, and $(s, c, S)_F$ policies with average cost improvements of $0.37\%$, $0.37\%$, $1.81\%$, $1.41\%$, and $1.67\%$; while the $(Q,S,T)$ policy performs slightly better than the $(R, S)$ policy with an average cost improvement of $0.19\%$.  

\begin{table*}[!ht]
  \centering
 \resizebox{0.95\linewidth}{!}{
    \begin{tabular}{rrr|r|rrrrrr}
    \hline
    \multicolumn{1}{r}{\multirow{2}[4]{*}{$K$}} & \multicolumn{1}{r}{\multirow{2}[4]{*}{$b$}} & \multicolumn{1}{r|}{\multirow{2}[4]{*}{$h$}} & \multicolumn{1}{c|}{\multirow{2}[4]{*}{$(R, S)$}} & \multicolumn{6}{c}{Average cost improvement $\Delta\%$} \\
\cline{5-10}          &       & \multicolumn{1}{r|}{} &       & $(Q, S, T)$ & $Q(s, S)$ & $P(s, S)$ & $(Q, S)$ & $MP$ & $(s, c, S)_F$ \\
    \hline
    \multicolumn{1}{r}{20} & \multicolumn{1}{r}{10} & \multicolumn{1}{r|}{2} & \multicolumn{1}{r|}{772.25} & \textcolor[rgb]{ 1,  0,  0}{\textbf{-0.03}} & 0.48  & 0.76  & 8.30  & 1.79  & 1.80  \\
    \multicolumn{1}{r}{50} & \multicolumn{1}{r}{10} & \multicolumn{1}{r|}{2} & \multicolumn{1}{r|}{813.94} & \textcolor[rgb]{ 1,  0,  0}{\textbf{-0.48}} & 0.12  & 0.62  & 0.47  & 1.64  & 1.74  \\
    \multicolumn{1}{r}{100} & \multicolumn{1}{r}{10} & \multicolumn{1}{r|}{2} & \multicolumn{1}{r|}{861.05} & 0.23  & 0.70  & 1.17  & 3.68  & 2.20  & 2.38  \\
    \multicolumn{1}{r}{200} & \multicolumn{1}{r}{10} & \multicolumn{1}{r|}{2} & \multicolumn{1}{r|}{932.86} & 1.62  & 1.83  & 2.38  & 2.88  & 3.42  & 3.73  \\
    \multicolumn{1}{r}{500} & \multicolumn{1}{r}{10} & \multicolumn{1}{r|}{2} & \multicolumn{1}{r|}{1131.42} & 0.14  & 0.14  & 0.59  & 0.18  & 1.60  & 2.12  \\
    \multicolumn{1}{r}{20} & \multicolumn{1}{r}{10} & \multicolumn{1}{r|}{6} & \multicolumn{1}{r|}{1166.06} & 0.85  & 2.84  & 0.01  & 7.99  & 1.08  & 1.04  \\
    \multicolumn{1}{r}{50} & \multicolumn{1}{r}{10} & \multicolumn{1}{r|}{6} & \multicolumn{1}{r|}{1222.82} & \textcolor[rgb]{ 1,  0,  0}{\textbf{-0.15}} & 1.83  & 0.62  & 5.53  & 1.68  & 1.73  \\
    \multicolumn{1}{r}{100} & \multicolumn{1}{r}{10} & \multicolumn{1}{r|}{6} & \multicolumn{1}{r|}{1283.92} & 1.33  & 2.50  & 1.26  & 4.49  & 2.34  & 2.46  \\
    \multicolumn{1}{r}{200} & \multicolumn{1}{r}{10} & \multicolumn{1}{r|}{6} & \multicolumn{1}{r|}{1413.72} & 0.30  & 1.23  & 1.02  & 1.82  & 2.10  & 2.33  \\
    \multicolumn{1}{r}{500} & \multicolumn{1}{r}{10} & \multicolumn{1}{r|}{6} & \multicolumn{1}{r|}{1658.48} & 2.26  & 2.20  & 2.52  & 2.30  & 3.59  & 4.03  \\
    \multicolumn{1}{r}{50} & \multicolumn{1}{r}{10} & \multicolumn{1}{r|}{10} & \multicolumn{1}{r|}{1420.63} & 1.57  & 5.30  & \textcolor[rgb]{ 1,  0,  0}{\textbf{-0.03}} & 5.88  & 1.07  & 1.07  \\
    \multicolumn{1}{r}{100} & \multicolumn{1}{r}{10} & \multicolumn{1}{r|}{10} & \multicolumn{1}{r|}{1497.96} & 1.67  & 4.28  & 0.75  & 4.37  & 1.87  & 1.93  \\
    \multicolumn{1}{r}{200} & \multicolumn{1}{r}{10} & \multicolumn{1}{r|}{10} & \multicolumn{1}{r|}{1637.27} & 0.66  & 2.18  & 1.15  & 2.16  & 2.28  & 2.44  \\
    \multicolumn{1}{r}{500} & \multicolumn{1}{r}{10} & \multicolumn{1}{r|}{10} & \multicolumn{1}{r|}{1935.07} & 1.60  & 1.60  & 1.79  & 1.60  & 2.90  & 3.27  \\
    \multicolumn{1}{r}{100} & \multicolumn{1}{r}{50} & \multicolumn{1}{r|}{2} & \multicolumn{1}{r|}{1043.31} & \textcolor[rgb]{ 1,  0,  0}{\textbf{-1.95}} & \textcolor[rgb]{ 1,  0,  0}{\textbf{-0.79}} & \textcolor[rgb]{ 1,  0,  0}{\textbf{-0.23}} & 1.98  & 0.78  & 0.92  \\
    \multicolumn{1}{r}{200} & \multicolumn{1}{r}{50} & \multicolumn{1}{r|}{2} & \multicolumn{1}{r|}{1132.61} & \textcolor[rgb]{ 1,  0,  0}{\textbf{-1.29}} & \textcolor[rgb]{ 1,  0,  0}{\textbf{-0.48}} & 0.30  & 0.50  & 1.31  & 1.97  \\
    \multicolumn{1}{r}{500} & \multicolumn{1}{r}{50} & \multicolumn{1}{r|}{2} & \multicolumn{1}{r|}{1327.95} & 0.08  & 0.08  & 0.82  & 0.13  & 1.83  & 2.30  \\
    \multicolumn{1}{r}{100} & \multicolumn{1}{r}{50} & \multicolumn{1}{r|}{6} & \multicolumn{1}{r|}{1794.60} & \textcolor[rgb]{ 1,  0,  0}{\textbf{-1.37}} & \textcolor[rgb]{ 1,  0,  0}{\textbf{-2.65}} & \textcolor[rgb]{ 1,  0,  0}{\textbf{-2.09}} & 0.94  & \textcolor[rgb]{ 1,  0,  0}{\textbf{-1.09}} & \textcolor[rgb]{ 1,  0,  0}{\textbf{-0.97}} \\
    \multicolumn{1}{r}{200} & \multicolumn{1}{r}{50} & \multicolumn{1}{r|}{6} & \multicolumn{1}{r|}{1938.25} & \textcolor[rgb]{ 1,  0,  0}{\textbf{-0.27}} & \textcolor[rgb]{ 1,  0,  0}{\textbf{-1.56}} & \textcolor[rgb]{ 1,  0,  0}{\textbf{-0.89}} & \textcolor[rgb]{ 1,  0,  0}{\textbf{-0.05}} & 0.13  & 0.34  \\
    \multicolumn{1}{r}{500} & \multicolumn{1}{r}{50} & \multicolumn{1}{r|}{6} & \multicolumn{1}{r|}{2244.01} & \textcolor[rgb]{ 1,  0,  0}{\textbf{-0.27}} & \textcolor[rgb]{ 1,  0,  0}{\textbf{-0.27}} & 0.43  & \textcolor[rgb]{ 1,  0,  0}{\textbf{-0.26}} & 1.44  & 1.87  \\
    \multicolumn{1}{r}{200} & \multicolumn{1}{r}{50} & \multicolumn{1}{r|}{10} & \multicolumn{1}{r|}{2448.79} & \textcolor[rgb]{ 1,  0,  0}{\textbf{-3.83}} & \textcolor[rgb]{ 1,  0,  0}{\textbf{-2.11}} & \textcolor[rgb]{ 1,  0,  0}{\textbf{-1.55}} & \textcolor[rgb]{ 1,  0,  0}{\textbf{-0.75}} & \textcolor[rgb]{ 1,  0,  0}{\textbf{-0.53}} & \textcolor[rgb]{ 1,  0,  0}{\textbf{-0.34}} \\
    \multicolumn{1}{r}{500} & \multicolumn{1}{r}{50} & \multicolumn{1}{r|}{10} & \multicolumn{1}{r|}{2796.29} & 0.35  & 0.35  & 0.97  & 0.35  & 2.00  & 2.40  \\
    \multicolumn{1}{r}{200} & \multicolumn{1}{r}{100} & \multicolumn{1}{r|}{2} & \multicolumn{1}{r|}{1200.38} & \textcolor[rgb]{ 1,  0,  0}{\textbf{-1.61}} & \textcolor[rgb]{ 1,  0,  0}{\textbf{-0.94}} & \textcolor[rgb]{ 1,  0,  0}{\textbf{-0.11}} & \textcolor[rgb]{ 1,  0,  0}{\textbf{-0.01}} & 0.90  & 1.13  \\
    \multicolumn{1}{r}{500} & \multicolumn{1}{r}{100} & \multicolumn{1}{r|}{2} & \multicolumn{1}{r|}{1406.67} & \textcolor[rgb]{ 1,  0,  0}{\textbf{-0.76}} & \textcolor[rgb]{ 1,  0,  0}{\textbf{-0.83}} & 0.16  & 0.16  & 1.17  & 1.60  \\
    \multicolumn{1}{r}{200} & \multicolumn{1}{r}{100} & \multicolumn{1}{r|}{6} & \multicolumn{1}{r|}{2106.78} & 0.44  & \textcolor[rgb]{ 1,  0,  0}{\textbf{-1.23}} & \textcolor[rgb]{ 1,  0,  0}{\textbf{-0.48}} & 0.94  & 0.54  & 0.73  \\
    \multicolumn{1}{r}{500} & \multicolumn{1}{r}{100} & \multicolumn{1}{r|}{6} & \multicolumn{1}{r|}{2449.51} & \textcolor[rgb]{ 1,  0,  0}{\textbf{-0.88}} & \textcolor[rgb]{ 1,  0,  0}{\textbf{-0.88}} & \textcolor[rgb]{ 1,  0,  0}{\textbf{-0.07}} & \textcolor[rgb]{ 1,  0,  0}{\textbf{-0.07}} & 0.94  & 1.33  \\
    \multicolumn{1}{r}{200} & \multicolumn{1}{r}{100} & \multicolumn{1}{r|}{10} & \multicolumn{1}{r|}{2728.08} & \textcolor[rgb]{ 1,  0,  0}{\textbf{-3.41}} & \textcolor[rgb]{ 1,  0,  0}{\textbf{-1.90}} & \textcolor[rgb]{ 1,  0,  0}{\textbf{-1.29}} & \textcolor[rgb]{ 1,  0,  0}{\textbf{-0.49}} & \textcolor[rgb]{ 1,  0,  0}{\textbf{-0.27}} & \textcolor[rgb]{ 1,  0,  0}{\textbf{-0.10}} \\
    \multicolumn{1}{r}{500} & \multicolumn{1}{r}{100} & \multicolumn{1}{r|}{10} & \multicolumn{1}{r|}{3108.05} & 0.22  & 0.22  & 0.94  & 0.94  & 1.96  & 2.33  \\
    \multicolumn{1}{r}{500} & \multicolumn{1}{r}{200} & \multicolumn{1}{r|}{2} & \multicolumn{1}{r|}{1470.29} & \textcolor[rgb]{ 1,  0,  0}{\textbf{-0.90}} & \textcolor[rgb]{ 1,  0,  0}{\textbf{-0.90}} & 0.05  & 0.05  & 1.05  & 1.45  \\
    \multicolumn{1}{r}{500} & \multicolumn{1}{r}{200} & \multicolumn{1}{r|}{6} & \multicolumn{1}{r|}{2620.77} & \textcolor[rgb]{ 1,  0,  0}{\textbf{-0.91}} & \textcolor[rgb]{ 1,  0,  0}{\textbf{-0.91}} & 0.08  & 0.08  & 1.09  & 1.45  \\
    \multicolumn{1}{r}{500} & \multicolumn{1}{r}{200} & \multicolumn{1}{r|}{10} & \multicolumn{1}{r|}{3421.28} & \textcolor[rgb]{ 1,  0,  0}{\textbf{-0.94}} & \textcolor[rgb]{ 1,  0,  0}{\textbf{-0.94}} & \textcolor[rgb]{ 1,  0,  0}{\textbf{-0.04}} & \textcolor[rgb]{ 1,  0,  0}{\textbf{-0.04}} & 0.97  & 1.30  \\
    \hline
    \multicolumn{4}{l|}{Average cost improvement $\Delta\%$} & \textcolor[rgb]{ 1,  0,  0}{\textbf{-0.19}} & 0.37  & 0.37  & 1.81  & 1.41  & 1.67  \\
\hline
\end{tabular}%
  }
\caption{Computational results on the data set of \cite{viswanathan1997}}
  \label{dataset2}%
\end{table*}%

Even though the $(R, S)$ policy does not fully dominate alternative policies, it presents a key advantage: {\em in contrast to all other policies in the literature, it is able to tackle stationary as well as nonstationary demand.}

\section{Conclusion}\label{conclusion}
In this paper, we presented a mathematical programming approach for controlling the multi-item inventory system with joint replenishment under the $(R, S)$ policy. We first presented an MILP-based model for approximating optimal $(R, S)$ policies, which is built upon the piecewise-linear approximation technique proposed by \citep{rossietal2014}.  We further demonstrated that the MILP model can be used to approximate the $(\sigma, \vec{S})$ policy.

We conducted an extensive computational study comprising $40$ instances. We first evaluated our approach on the data set of \cite{atkinsandIyogun1988}. This evaluation demonstrates that the $(R, S)$ policy fully dominates other competing policies in the literature in 2 out of $9$ test instances considered. The $(R, S)$ policy performs better than $Q(s, S)$, $(Q, S)$, $MP$, $(s, c, S)_M$, and $(s, c, S)_F$ policies with  average improvements of $0.07\%$, $1.74\%$, $0.89\%$, $2.84\%$, and $7.02\%$, respectively; however, $(Q, S, T)$ and $P(s, S)$ policies perform slightly better than the $(R, S)$ policy with average improvements of $0.09\%$ and $0.14\%$.  Computational experiments on the data set of \cite{viswanathan1997} indicate that $(R, S)$ is the best policy in $13$ out of $31$ test instances. The $(R, S)$ policy performs better than $Q(s, S)$,  $P(s, S)$, $(Q, S)$, $MP$, and $(s, c, S)_F$ policies with average cost improvements of $0.37\%$, $0.37\%$, $1.84\%$, $1.41\%$ and $1.67\%$; while the $(Q,S,T)$ policy performs slightly better than the $(R, S)$ policy with an average cost improvement $0.19\%$. Most importantly, the $(R, S)$ policy comes with the additional advantage of being able to tackle stationary and nonstationary demand. Future research may focus on investigating the cost performance of $(R, S)$ policy in a rolling horizon setting. 





\bibliography{sample}







\appendix
\section{MILP model for approximating $(R, S)$ policies with fixed lead time}\label{leadtime}

This section demonstrates that the MINLP model in Fig. \ref{MINLPmodel} can be extended to compute near-optimal $(R, S)$ policy parameters for nonstationary JRPs with fixed lead time. Let $L^n$ denote the lead time of item $n$, $n=1, \ldots, N$. We next separate our discussions into two parts. 

The first part involves periods $1, \dots, L^n$, $n=1, \ldots, N$, where no order is received. We assume that there is no outstanding order at the beginning of the planning horizon, and the system is forced to issue an order in period $1$, then the inventory level $I_t$ must equal to the initial inventory level of item $n$ at the beginning of the planning horizon, minus the demand convolution over periods $1, \ldots, t$, i.e., $I_t^n=I_0^n -d_{1,t}^n$, where $d_{1,t}^n$ is the demand convolution of item $n$ over periods $1, \ldots, t$, i.e., $d_{1,t}^n=d_1^n+\ldots +d_t^n$. We rewrite the expected back-orders and excess on-hand stocks using the first order loss function and its complementary function, $\mathcal{L}(I_0^n, d_{1,t}^n)$ and $\hat{\mathcal{L}}(I_0^n, d_{1,t}^n)$. 

Additionally, since no order of item $n$ is received before $L^n$, the expected inventory level of item $n$ at the end of period $t$ is equal to the expected inventory level at the end of period $t-1$, minus expected demand in period $t$, $t=1, \ldots, L^n$, 
\begin{align}
&\tilde{I}_t^n+\tilde{d}_t^n-\tilde{I}_{t-1}^n=0, &t={1, \ldots, L^n}.
\end{align}

The second part involves periods $1+L^n, \ldots, T$, $n=1, \ldots, N$. Consider a single cycle of item $n$ over periods $i,\ldots, j$, in which a single order is received at the beginning of period $i$, and the next order will be received at the beginning of period $j+1$. Since the lead time of item $n$ is $L^n$, the order that arrives in period $i$ must be issued in period $i-L^n$ with order-up-to-position $S_{i-L^n}^n$. Thus, $I_t^n$, $t=i, \ldots, j$, must equal to the order-up-to-position $S_{i-L^n}^n$, minus the demand convolution over periods $i-L^n, \ldots, t$, i.e. $I_t^n=S_{i-L^n}^n-d_{i-L^n,t}^n$. 

We introduce a binary variable $P_{jt}^n$ which is set to one if the most recent order received before period $t$ arrived in period $j$, where $j\leq t$, $j=1+L^n, \ldots, t$, $t=1+L^n, \ldots, T$, and $n=1, \ldots, N$; and we introduce the following constraints, $t={1+L^n, \ldots, T}$, $n=1, \ldots, N$,
\begin{align}
&\sum_{j=1+L^n}^tP_{jt}^n=1, & \label{c-1}\\
&P_{j,t}^n \geq y_{j-L^n}^n - \sum_{k=j-L^n+1}^{t-L^n}y_k^n, & j={1+L^n, \ldots, t} \label{c-2}.
\end{align}
Constraints (\ref{c-1}) indicate that the most recent order received before period t arrived in period $j$. Constraints (\ref{c-2}) identify uniquely the period in which the most recent order received before period t has been received. Therefore, the inventory level $I_t^n=\sum_{j=1+L^n}^t(S_{j-L^n}^n - d_{j-L^n,t}^n)P_{jt}^n$, where $t=1+L^n, \ldots, T$, and $S_{j-L^n}^n$ represents the order-up-to-position of item $n$ in period $j-L^n$. We write the back-orders and excess inventory  as the first order loss function and its complementary, $\sum_{j=1+L^n}^t\mathcal{L}(S_{j-L^n}^n, d_{j-L^n,t}^n)P_{jt}^n$ and $\sum_{j=1+L^n}^t\hat{\mathcal{L}}(S_{j-L^n}^n, d_{j-L^n,t}^n)P_{jt}^n$. 


In addition, constraints (\ref{MINLP-2})-(\ref{MINLP-3}) in Fig. \ref{MINLPmodel} can be reformulated as follows,
\begin{align}
&y^n_{t-L^n}=0 \rightarrow \tilde{I}_{t}^n+\tilde{d}_{t}^n-\tilde{I}_{t-1}^n=0, &t={1+L^n, \ldots, T}, \\
&\tilde{I}_t^n+\tilde{d}_t^n-\tilde{I}_{t-1}^n\geq 0, &t={1+L^n, \ldots, T}.
\end{align}

We now present the MILP model for approximating $(R, S)$ policies with fixed lead time in Fig. \ref{MILP}. 
\begin{figure}[!ht]
\tiny
\begin{equation}
\min \sum_{t=1}^T\Big(K\cdot \delta_{t}+\sum_{n=1}^N\big(k^n\cdot y_t^n+h^n\tilde{H}_{t}^n+b^n\tilde{B}_{t}^n\big)\Big) \label{MILP-0}
\end{equation}
Subject to, $n = 1, \ldots, N$
\begin{align}
&\delta_t\geq y_t^n & t=1, \ldots, T 
\label{MILP-1}\\
&y_1^n=1 & 
\label{MILP-2}\\
&\tilde{I}_t^n+\tilde{d}_t^n-\tilde{I}_{t-1}^n=0 &t={1, \ldots, L^n} 
\label{MILP-3}\\
&B_t^n \geq -\tilde{I}_t^n+\sum_{k=1}^ip_kI_0^n-\sum_{k=1}^ip_kE[d_{1,t}^n|\Omega_i], & \begin{array}{l} t=1,\ldots, L^n\\ i=1, \ldots, W \end{array}
\label{MILP-6}\\
&H_t^n \geq \sum_{k=1}^ip_kI_0^n-\sum_{k=1}^ip_kE[d_{1,t}^n|\Omega_i], &\begin{array}{l} t=1,\ldots, L^n \\ i=1, \ldots, W \end{array}
\label{MILP-7}\\
&y_t^n=0 & t=T-L^n, \ldots, T \label{MILP4-1}\\
&\tilde{I}_t^n+\tilde{d}_t^n-\tilde{I}_{t-1}^n\geq 0 &t={1+L^n, \ldots, T}
\label{MILP-4}\\
&y^n_{t-L^n}=0 \rightarrow \tilde{I}_{t}^n+\tilde{d}_{t}^n-\tilde{I}_{t-1}^n=0 &t={1+L^n, \ldots, T}  \label{MILP-5}\\
&\sum_{j=1+L^n}^tP_{jt}^n=1 & t={1+L^n, \ldots, T} 
\label{MILP-8}\\
&P_{j,t}^n \geq y_{j-L^n}^n - \sum_{k=j-L^n+1}^{t}y_k^n & \begin{array}{l} t=1+L^n, \ldots, T \\ j=1, \ldots, t \end{array}\label{MILP-9}\\
&B_t^n \geq -\tilde{I}_t^n+(\tilde{I}_t^n+\sum_{j=1+L^n}^t\tilde{d}_{j-L^n,t}^nP_{jt}^n)\sum_{k=1}^ip_k-\sum_{j=1+L^n}^t\sum_{k=1}^ip_kE[d_{j-L^n,t}^n|\Omega_i]P_{jt}^n & \begin{array}{l} t=1+L^n, \ldots, T\\ i=1, \ldots, W\end{array}
\label{MILP-10}\\
&H_t^n \geq (\tilde{I}_t^n+\sum_{j=1+L^n}^t\tilde{d}_{j-L^n,t}^nP_{jt})\sum_{k=1}^ip_k-\sum_{j=1+L^n}^t\sum_{k=1}^ip_kE[d_{j-L^n,t}^n|\Omega_i]P_{jt}^n &\begin{array}{l} t=1+L^n, \ldots, T\\ i=1, \ldots, W \end{array}
\label{MILP-11}\\
&\delta_t=\{0, 1\} & t={1, \ldots, T} 
\label{MILP-12}\\
&y_t^n =\{0, 1\} &t={1, \ldots, T} 
\label{MILP-13}\\
&P_{jt}^n=\{0, 1\} & \begin{array}{l} t={1+L^n, \ldots, T} \\j={1+L^n, \ldots, t} \end{array}
\label{MILP-14}
\end{align}
\caption{MILP model for approximating $(R, S)$ policies}
\label{MILP}
\end{figure}

The objective function (\ref{MILP-0}) minimises the expected group fixed ordering costs, item-specific fixed ordering costs, penalty costs, and holding costs of $N$-item over the $T$-period planning horizon. 
Constraints (\ref{MILP-1}) imply  an individual item can only be included in a group replenishment if that replenishment is made.
Constraints (\ref{MILP-2}) - (\ref{MILP-3}) assume that the first order is issued at the beginning of period $1$, and there is no outstanding replenishment at the beginning of the planning horizon. 
Constraints (\ref{MILP-6}) - (\ref{MILP-7}) represent the expected back-orders and on-hand stocks of item $n$ over periods $1, \ldots, L^n$.
Constraints (\ref{MILP4-1}) state all orders are received by the end of the planning horizon.
Constraints (\ref{MILP-4}) - (\ref{MILP-5}) are inventory balance constraints. 
Constraints (\ref{MILP-8}) - (\ref{MILP-9}) ensure the most recent replenishment that has arrived before period $t$ was received in period $j$.
Constraints (\ref{MILP-10}) - (\ref{MILP-11}) represent the expected back-orders  and  on-hand stocks of item $n$ over periods $1+L^n, \ldots, T$.
Constraints (\ref{MILP-12}) - (\ref{MILP-14}) indicate domains of binary variables $\delta_t^n$, $y_t^n$, and $P_{jt}^n$.

\section{Shortest path reformulation for approximating stationary stochastic $(R, S)$ policies} \label{shortespathapproximation}
In this section we present an efficient shortest path reformulation for computing stationary $(R, S)$ policies.

Consider a network $\mathcal{G}=(\mathcal{N}, \mathcal{A})$ with nodes $\mathcal{N} = \{1, \ldots, T\}$ representing time periods, and arc $(i, j)$ between each pair of $(i, j)$ representing a possible decision to issue an order  in period $i$ to satisfy demands in periods $i, \ldots, j$. Assigning a cost to this arc,  solving the optimisation problem in Fig. \ref{MINLPmodel} is equivalent to finding the shortest path between nodes $1$ and $T$ in the network $\mathcal{G}$. 
In the rest of this section, we first present how to compute the cost of each arc, and then present the shortest path reformulation. 

Consider a replenishment cycle $i, \ldots, j$, where the only order is issued in period $i$ with order-up-to-position $S_{ij}^n$, and the next order is issued in period $j+1$, for $i=1, \ldots, T$, $j=i, \ldots, T$, $n=1, \ldots, N$. We assume $d_t^n$ follows a Poisson distribution with rate $\lambda^n$. Then, $S_{ij}^n$ can be calculated by \cite{askin1981},
\begin{align}
   & \sum_{t=i}^j G_{d_{i,t+L^n}^n}(S_{i,j}^n) = \frac{(j-i+1)\cdot b^n}{h^n+b^n}.\label{calculateS}
\end{align}
Note that the order-up-to-position $S_{i,j}^n$ actually accounts for demand variances over periods $i, \ldots, j+L^n$, which is reflected on the cumulative distribution function $G_{d_{1,t+L^n}^n}(\cdot)$ on the left-hand-side of Eq. (\ref{calculateS}). 

Since the demand of item $n$ follows a Poisson distribution with rate $\lambda^n$, we could approximate the cost of the replenishment cycle $i, \ldots, j$ by that of the cycle $i+L^n, \ldots, j+L^n$ as shown in Fig. \ref{stationaryreplenishment}. As a result, the cycle cost $c_{ij}^n$ can be calculated as follows,
\begin{align}
    &c_{ij}^n=k^n
    +h^n \sum_{t=i}^j\hat{\mathcal{L}}(S^n_{i,j}-L^n\lambda^n, d_{it})
    +b^n\sum_{t=i}^j\mathcal{L}(S^n_{i,j}-L^n\lambda^n, d_{it}).
\end{align}
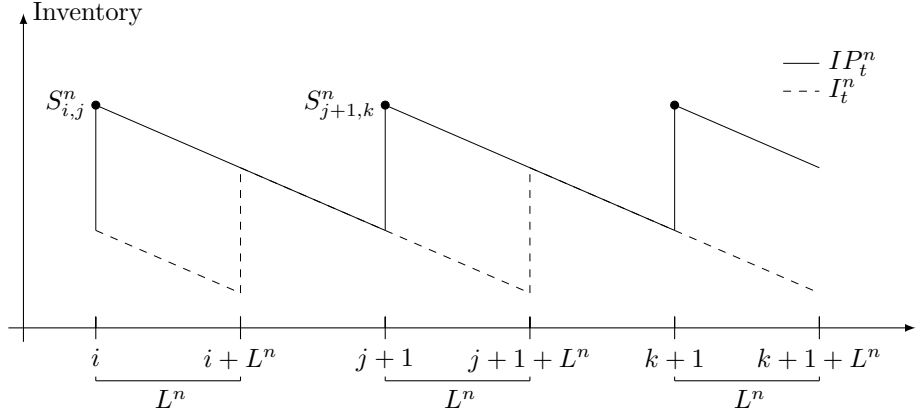
\begin{figure}[!ht]
\centering
\begin{tikzpicture}[x=0.0480952380952381cm, y=0.018518518518518517cm]
 \draw [-latex] ([xshift=-2mm] 0.0,0) -- ([xshift=115mm] 7.5,0);
 \draw (20.0,0) -- +(0mm,1.5mm) -- +(0mm,-1.5mm) node[below] {$i$};
 \draw (60,0) -- +(0mm,1.5mm) -- +(0mm,-1.5mm) node[below] {$i+L^n$};
 \draw (100,0) -- +(0mm,1.5mm) -- +(0mm,-1.5mm) node[below] {$j+1$};
 \draw (140,0) -- +(0mm,1.5mm) -- +(0mm,-1.5mm) node[below] {$j+1+L^n$};
 \draw (180,0) -- +(0mm,1.5mm) -- +(0mm,-1.5mm) node[below] {$k+1$};
 \draw (220,0) -- +(0mm,1.5mm) -- +(0mm,-1.5mm) node[below] {$k+1+L^n$};
 \draw [-latex] ([yshift=-0mm] 0,-10.0) -- ([yshift=40mm] 0, 10.0) node[right] {Inventory};

 \filldraw[black] (20,160) circle (1.5pt) node[anchor=east] {$S_{i,j}^n$};
 \filldraw[black] (100,160) circle (1.5pt) node[anchor=east] {$S_{j+1,k}^n$};
 \filldraw[black] (180,160) circle (1.5pt) node[anchor=east] {};
 
 \draw [color=black, mark= , style=solid](20,-35)--(20,-38);
\draw [color=black, mark= , style=solid](20,-38) node[below,near start abs] {$L^n$}--(60,-38);
\draw [color=black, mark= , style=solid](60,-35)--(60,-38);

 \draw [color=black, mark= , style=solid](100,-35)--(100,-38);
\draw [color=black, mark= , style=solid](100,-38) node[below,near start abs] {$L^n$}--(140,-38);
\draw [color=black, mark= , style=solid](140,-35)--(140,-38);

  \draw [color=black, mark= , style=solid](180,-35)--(180,-38);
 \draw [color=black, mark= , style=solid](180,-38) node[below,near start abs] {$L^n$}--(220,-38);
 \draw [color=black, mark= , style=solid](220,-35)--(220,-38);

 \draw[line width=0.1mm, color=black, mark= , style=solid](210,190)--(220,190)node[pos=1,right]{\textcolor{black}{$IP_t^n$}};

 \draw [line width=0.1mm, color=black, mark=,style=solid](20,70)--(20,160);
 \draw [line width=0.1mm, color=black, mark= , style=solid] (20,160)--(100,70);
 \draw [line width=0.1mm, color=black, mark= , style=solid](100,70)--(100,160);
 \draw [line width=0.1mm, color=black, mark= , style=solid] (100,160)--(180,70);
 \draw [line width=0.1mm, color=black, mark= , style=solid](180,70)--(180,160);
 \draw [line width=0.1mm, color=black, mark= , style=solid] (180,160)--(220,115);
 
  \draw[line width=0.1mm, color=black, mark= , style=dashed](210,170)--(220,170)node[pos=1,right]{\textcolor{black}{$I_t^n$}};
 \draw [line width=0.1mm, color=black, mark= , style=dashed] (20,70)--(60,25);
 \draw [line width=0.1mm, color=black, mark= , style=dashed] (60,25)--(60,115);
  \draw [line width=0.1mm, color=black, mark= , style=dashed] (60,115)--(100,70);
 \draw [line width=0.1mm, color=black, mark= , style=dashed] (100,70)--(140,25);
  \draw [line width=0.1mm, color=black, mark= , style=dashed] (140,25)--(140,115);
  \draw [line width=0.1mm, color=black, mark= , style=dashed] (140,115)--(180,70);
  \draw [line width=0.1mm, color=black, mark= , style=dashed] (180,70)--(220,25);
\end{tikzpicture}
\caption{Expected inventory curve under $(R, S)$ policy.}
\label{stationaryreplenishment}
\end{figure}

At the beginning of the planning horizon, the initial inventory level is $I_0^n$. We check the cost of not issuing an order in period $1$, $\bar{c}_{1j}^n$, and update $c_{1j}^n$ with $\bar{c}_{1j}^n$ if $\bar{c}_{1j}^n \leq c_{1j}^n$, for $j=1, \ldots, T$.
\begin{align}
    &\bar{c}^n_{1j}=
    h^n\cdot \sum_{t=1}^j\hat{\mathcal{L}}(I_0^n, d_{1t})
    +b^n\cdot \sum_{t=1}^j\mathcal{L}(I_0^n, d_{1t}).
\end{align}
Additionally, we introduce an auxiliary binary variable $P_j^n$, which is equal to $1$ if an order is placed in period $1$ to satisfy demands in cycle $1, \ldots, j$, otherwise 0.

We now present the shortest path reformulation in Fig. \ref{shortestpath}. Let binary variable $Y_{ij}^n$ equal to $1$ if an order is issued in period $i$ to cover demands in periods $i, \ldots, j$, otherwise 0. The objective is to find the optimal replenishment plan that minimising the expected group fixed order costs, item-specific fixed order costs, holding costs and penalty costs over periods $1, \ldots, T$ for items $1, \ldots, N$.
\begin{figure}[!ht]
\tiny
\begin{align}
& \min \sum_{i=1}^T K\cdot \delta_i + \sum_{n=1}^N\sum_{i=1}^T\sum_{j=i}^Tc^n_{ij}\cdot Y^n_{ij} \label{shortestpathobj}
\end{align}
subject to, $n=1, \ldots, N$,
\begin{align}
&\delta_1 \geq \sum_{j=1}^TY_{1j}^n\cdot P_j^n &\label{shortestpath-4}\\
&\delta_i \geq \sum_{j=i}^TY_{ij}^n & i=2, \ldots, T \label{shortestpath-5}\\
& \sum_{j=1}^TY_{1j}^n = 1 & \label{shortestpath-1}\\
& \sum_{j=i}^TY_{ij}^n - \sum_{k=1}^{i-1}Y_{ki}^n = 0 & i=2, \ldots, T-1\label{shortestpath-2}\\
&\sum_{i=1}^TY_{iT}^n = 1 &\label{shortestpath-3}
\end{align}
\caption{Shortest path formulation for approximating stationary stochastic $(R, S)$ policies}
\label{shortestpath}
\end{figure}

Recall that $P_j^n$ represents the item-specific first period replenishment decision, which is set to $1$ if an order is issued in period $1$, otherwise $0$. Therefore, Constraints (\ref{shortestpath-4}) guarantee the group fixed order cost in period $1$ is properly counted. Constraints (\ref{shortestpath-5}) ensure that the group fixed order cost is encountered whenever any item is replenished in period $2, \ldots, T$. Constraints (\ref{shortestpath-1}) ensure that there is no more than one outgoing arc from period $1$. Constraints (\ref{shortestpath-2}) are flow balance equations. Constraints (\ref{shortestpath-3}) guarantee that period $T$ is included in a replenishment cycle. By solving the shortest path reformulation in Fig. \ref{shortestpath}, the group order decision $\delta_t^n$ and item-specific order decision $y_t^n$ are obtained\footnote{This can be obtained by adding constraints $y_1^n = \sum_{j=1}^TY_{1j}^nP_j^n$ and $y_i^n=\sum_{j=2}^TY_{ij}^n$, $i=2, \ldots, T$, to Fig. \ref{shortestpath}.}, for $t=1, \ldots, T$, $n = 1, \ldots, N$.

\end{document}